\documentclass[a4paper,twoside,11pt,reqno]{amsart}

\usepackage[headings]{fullpage}
\usepackage{lmodern}
\usepackage{microtype}
\usepackage{amsmath}
\usepackage{amsthm}
\usepackage{amssymb}
\usepackage{xcolor}
\usepackage{hyperref}
\hypersetup{colorlinks=true,linkcolor=red,citecolor=black,urlcolor=black}
\usepackage{mathrsfs}
\usepackage{enumitem}
\usepackage{tikz-cd}

\theoremstyle{definition}
\newtheorem{defi}{Definition}[section]

\newtheorem{rem}[defi]{Remark}
\newtheorem{defilem}[defi]{Definition and Lemma}
\theoremstyle{plain}
\newtheorem*{thm*}{Theorem}
\newtheorem{thm}[defi]{Theorem}
\newtheorem{lem}[defi]{Lemma}
\newtheorem{cor}[defi]{Corollary}
\newtheorem{prop}[defi]{Proposition}
\newtheorem*{conj*}{Conjecture}

\newcommand{\ab}{\mathrm{ab}}
\renewcommand{\ss}{\mathrm{ss}}

\newcommand{\pr}{\mathrm{pr}}

\newcommand{\cusp}{\mathrm{cusp}}
\newcommand{\rat}{\mathrm{rat}}

\renewcommand{\P}{\mathbb{P}}

\newcommand{\N}{\mathbb{N}}
\newcommand{\Q}{\mathbb{Q}}
\newcommand{\Z}{\mathbb{Z}}
\newcommand{\F}{\mathbb{F}}
\newcommand{\A}{\mathbb{A}}
\newcommand{\id}{\mathrm{id}}
\renewcommand{\H}{\mathrm{H}}
\renewcommand{\o}{\mathcal{O}}

\newcommand{\lisim}{\overset{\sim}{\longrightarrow}}
\DeclareMathOperator{\GL}{GL}

\DeclareMathOperator{\SL}{SL}
\DeclareMathOperator{\spec}{Spec}
\DeclareMathOperator{\im}{im}

\DeclareMathOperator{\Hom}{Hom}
\DeclareMathOperator{\Gal}{Gal}
\DeclareMathOperator{\Aut}{Aut}

\DeclareMathOperator*{\colim}{colim}
\newcommand{\Ell}{\mathrm{Ell}}
\newcommand{\iso}{\simeq}
\newcommand{\bp}{\begin{proof}}
\newcommand{\ep}{\end{proof}}
\newcommand{\be}{\begin{enumerate}}
\newcommand{\ee}{\end{enumerate}}
\renewcommand{\epsilon}{\varepsilon}
\renewcommand{\phi}{\varphi}
\renewcommand{\l}{\ell}
\newcommand\restr[2]{{
		\left.\kern-\nulldelimiterspace 
		#1 
		\vphantom{\big|} 
		\right|_{#2} 
}}

\newcommand{\U}{\mathcal{U}}

\newcommand{\X}{\mathcal{X}}
\newcommand{\V}{\mathcal{V}}
\newcommand{\Y}{\mathcal{Y}}
\newcommand{\Tate}{\mathrm{Tate}}

\begin{document}
\title{Fermat's Last Theorem for Selmer sections}
\author{Benjamin Steklov}
\address{Benjamin Steklov, Institut f\"ur Mathematik, Goethe--Universit\"at Frankfurt, Robert-Mayer-Str.~6--8, 60325~Frankfurt am Main, Germany}
\email{steklov@math.uni-frankfurt.de}
\date{\today} 

		\begin{abstract}
We prove, in the context of the section conjecture, that every Selmer section over $\Q$ of the affine Fermat curve with prime exponent $\l\geq 7$ is cuspidal.
	\end{abstract}
	\maketitle
	\tableofcontents
\section{Introduction}
Let $V/k$ be a smooth, geometrically connected, hyperbolic curve over a field of characteristic zero $k$ with algebraic closure $\bar{k}$. Associated to $V/k$, we have a short exact fibration sequence of \'etale fundamental groups
\[1 \to \pi_1(\overline{V})\to \pi_1(V)\to \Gal_k \to 1\]
with $\overline{V}=V_{\bar{k}}$.
Let $\mathscr{S}_{\pi_1(V/k)}$ denote the set of sections of this sequence, considered up to $\pi_1(\overline{V})$-conjugacy. By the functoriality of $\pi_1$, a $k$-rational point $x\in V(k)$ determines a well-defined conjugacy class $s_x\in\mathscr{S}_{\pi_1(V/k)}$. This construction defines the \textit{non-abelian Kummer map}
\[\kappa: V(k)\longrightarrow \mathscr{S}_{\pi_1(V/k)},\; x \longmapsto s_x. \] 
Sections in the image of $\kappa$ are called \textit{rational}. When $V/k$ is not necessarily proper with smooth completion $V\subset X$ and boundary $Y=X\setminus V$, a $k$-rational cusp $y \in Y(k)$ gives rise to a split exact sequence 
\[1 \to \pi_1(\overline{V}_{\!y})\to \pi_1(V_y)\to \Gal_k\to 1,\]
where $V_y=V \times_X \spec(\o_{X,y}^h)$ is the \textit{scheme of nearby points}. The $k$-morphism $V_y\to V$ induces a map of sections $\mathscr{S}_{\pi_1(V_y/k)}  \to \mathscr{S}_{\pi_1(V/k)}$ and sections in the image for some $y \in Y(k)$ are called \textit{cuspidal}. In his letter sent to Faltings in 1983 \cite{grothendieckLetterFaltingsTranslation1997}, Grothendieck conjectured the following.
\newtheorem*{sectionconj}{Section conjecture}
\begin{sectionconj}[Grothendieck]
	Let $V/k$ be a smooth, geometrically connected, hyperbolic curve over a number field $k$. Then the map
		\[\kappa: V(k)\sqcup \bigsqcup_{y \in Y(k)}\mathscr{S}_{\pi_1(V_y/k)}\longrightarrow \mathscr{S}_{\pi_1(V/k)}\]
		is bijective.
\end{sectionconj}
 In other words, over a number field $k$, every section $s \in \mathscr{S}_{\pi_1(V/k)}$ is either explained by a unique rational point $x \in V(k)$ or belongs to a \textit{packet} of cuspidal sections  $\mathscr{S}_{\pi_1(V_y/k)}$ based at a unique $k$-rational cusp $y\in Y(k)$.
 
 \vspace{\baselineskip}
In this generality the section conjecture is still widely open. Injectivity was already known to Grothendieck, while the surjectivity of $\kappa$ remains the hard part. Some progress has been made on \textit{Selmer sections}, defined as follows. Let $k$ be a number field and $v$ a place of $k$ with completion $k_v$. A section $s \in \mathscr{S}_{\pi_1(V/k)}$ can be base-changed to a section $s_{k_v} \in \mathscr{S}_{\pi_1(V_{k_v}/k_v)}$. We call $s \in \mathscr{S}_{\pi_1(V/k)}$ a \textit{Selmer} section, if for all places $v$ the base-change $s_{k_v}$ is rational or cuspidal and we denote the corresponding subset by $ \mathscr{S}^\mathrm{Sel}_{\pi_1(V/k)} \subset\mathscr{S}_{\pi_1(V/k)}$.  In particular, the section conjecture predicts the following.
\newtheorem*{selsectionconj}{Selmer section conjecture}
 \begin{selsectionconj}
 	Let $V/k$ be a smooth, geometrically connected, hyperbolic curve over a number field $k$. Then the map
 	\[\kappa: V(k)\sqcup \bigsqcup_{y \in Y(k)}\mathscr{S}_{\pi_1(V_y/k)}\longrightarrow \mathscr{S}^\mathrm{Sel}_{\pi_1(V/k)}\]
 	is bijective.
 \end{selsectionconj} 
Stoll's results on the finite descent obstruction \cite[Corollary 8.8.]{stollFiniteDescentObstructions2007} imply the Selmer section conjecture for the modular curves $X_0(N)$ and $X_1(N)$ over $\Q$ whenever these curves have positive genus, see Proposition \ref{prop:finitedescentmoduli}. Using these results, we establish a new instance of the Selmer section conjecture, namely for the affine Fermat curves.
 \newtheorem*{theoremA}{Theorem A}
\begin{theoremA}[see Theorem \ref{thm:mainresult}]
Let $\l\geq 7$ be a prime and $V_\l/\Q$ the affine curve given by 
\[X^\l +Y^\l +Z^\l =0 \quad \text{and} \quad XYZ \neq 0\]
in the projective plane. Every Selmer section $s \in  \mathscr{S}_{\pi_1(V_\l/\Q)}$ is cuspidal.
\end{theoremA}
 Although $V_\ell$ has several $\Q$-rational cusps, Fermat's Last Theorem, proven by Wiles in 1994 \cite{wilesModularEllipticCurves1995}, asserts that $V_\ell$ has no rational points over $\Q$. Consequently, Theorem A establishes the Selmer section conjecture for the affine Fermat curves with exponent $\l\geq 7$.
 
Theorem A does not provide a new proof of Fermat's Last Theorem. In fact, its proof relies crucially on Serre's modularity conjecture, proven by Khare and Wintenberger \cite{khareSerresModularityConjecture2009}. Rather, it shows that techniques based on modularity can also apply to Selmer sections. We now explain the proof of Theorem A and the reason for the assumption $\l\geq 7$.
\subsection{Method of proof}

Our proof is based on Serre’s approach to Fermat’s Last Theorem.  In 1987, Serre showed that his conjecture on the modularity of residual $\Gal_{\Q}$-representations implies Fermat's Last Theorem for $\l\geq 5$ \cite{serreRepresentationsModulairesDegre1987}. To a rational point of the Fermat curve, one associates a corresponding \textit{Frey curve}, which is an elliptic curve over $\Q$ with affine equation 
\[E\colon y^2=x(x-a^\l)(x+b^\l),\]
where $a,b$ are integers satisfying $a^\l+b^\l=c^\l$ for some integer $c$. After several reduction steps involving the congruence properties of $a$ and $b$ modulo $8$ and $32$, the elliptic curve $E$ exhibits unusual reduction behavior. Applying Serre's modularity conjecture to the Galois representation $\overline{\rho}_{E,\l}$ associated to the $\l$-torsion points of $E$ yields a contradiction.

\vspace{\baselineskip}

Since Serre's modularity conjecture is now a theorem, proven by Khare and Wintenberger \cite{khareSerresModularityConjecture2009} in 2008, we adapt this strategy to Selmer sections. In place of a rational point, we attach a residual Galois representation to a section $s \in \mathscr{S}_{\pi_1(V_\l/\Q)}$. To this end, it is necessary to consider families of elliptic curves. More generally, given a hyperbolic curve $V/\Q$ and a family of elliptic curves $\mathcal{E}/V$, the $\l$-torsion of $\mathcal{E}$ yields a $\pi_1$-representation 
\[  \overline{\rho}_{\mathcal{E},\l}\colon \pi_1(V)\to \GL_2(\F_\l).\]
For a section $s \in \mathscr{S}_{\pi_1(V/\Q)}$, we obtain the Galois representation
\[\overline{\rho}_{s,\l}:=\overline{\rho}_{\mathcal{E},\l}\circ s \colon \Gal_{\Q} \to \GL_2(\F_\l).\]
The method of attaching Galois representations arising from families of elliptic curves or more general local systems to sections has been successfully used to prove results about Selmer sections, see \cite{stixBirationalSectionConjecture2015} and more recently \cite{bettsGaloisSectionsPadic2025}. In the case of the Fermat curve, we exploit the fact that $V_\l$ is a finite \'etale cover of the projective line minus three points $\P^1_\Q \setminus \{0,1,\infty\}$. The affine hyperbolic curve  $\P^1_\Q \setminus \{0,1,\infty\}$ carries the \textit{Legendre family} of elliptic curves, defined by the simple affine equation 
\[y^2 =x(x-1)(x-\lambda) \quad \lambda \in \P^1_\Q \setminus \{0,1,\infty\}.\]
The relation to the Frey curve is explained in the next subsection. Given a Selmer section $s \in \mathscr{S}_{\pi_1(V_\l/\Q)}$, we study the residual Galois representation $\overline{\rho}_{s,\l}$ associated with the pullback of the Legendre family to $V_\l$. The assumption that $s$ is a Selmer section is used crucially to control the ramification of $\overline{\rho}_{s,\l}$ and relate $s$ to the \textit{finite descent obstruction}. 

\vspace{\baselineskip}

The proof splits into two cases, depending on whether $\overline{\rho}_{s,\l}$ is irreducible or not. When $\overline{\rho}_{s,\l}$ is irreducible, we use Serre's modularity conjecture to derive a contradiction. In the case that  $\overline{\rho}_{s,\l}$ is reducible, we rely on Mazur's results on modular curves \cite{mazurModularCurvesEisenstein1977}, as well as Stoll's subsequent work on the finite descent obstruction \cite{stollFiniteDescentObstructions2007} to conclude that $s$ is cuspidal. Our method breaks down for $\l=5$, because the modular curve $X_0(10)$ has genus zero.
 \subsection{The use of DM-stacks}
Working with the Legendre family rather than a Frey curve poses several technical obstructions. In general, Frey curves are not part of the Legendre family, but are quadratic twists thereof.  
Indeed, over $\Q$, an elliptic curve of the form
\[
y^2 = x(x-a)(x+b)
\]
can be transformed into Legendre form
\[
y^2 = x(x-1)(x-\lambda)
\]
only after passing to a suitable quadratic extension of $\Q$. Consequently, instead of working with a fixed section $s \in \mathscr{S}_{\pi_1(V_\l/\Q)}$, we would at least like to consider the set of all ``quadratic twists'' of $s$. To make this possible, we push the method of using families of elliptic curves even further: In place of a fixed family, we use the entire moduli stack of elliptic curves $\Ell_{\Q}$. 

\vspace{\baselineskip}

Let $k$ be a field of characteristic zero. The moduli stack of elliptic curves $\Ell_{k}$ over $k$ is an example of a \textit{hyperbolic DM-curve} over $k$. By this we mean a one-dimensional DM-stack $\U$, which is (geometrically) connected, smooth and separated over $k$, admitting a finite \'etale cover by a hyperbolic curve. The section conjecture for DM-stacks $\U/k$ has been studied by Bresciani \cite{brescianiImplicationsGrothendiecksAnabelian2021} using the \textit{\'etale fundamental gerbe} $\Pi_{\U/k}$, with a focus on hyperbolic \textit{orbicurves}, i.e.\ hyperbolic DM-curves which admit a dense open subscheme. To accommodate the moduli stack $\Ell_k$, which is not an orbicurve, we formulate the section conjecture for arbitrary hyperbolic DM-curves. Building on Bresciani's work, we obtain the following generalization of \cite[Theorem~7.2(2)]{brescianiImplicationsGrothendiecksAnabelian2021}.
\newtheorem*{theoremB}{Theorem B}
\begin{theoremB}[see Theorem \ref{thm:abstractff}]
The section conjecture for all hyperbolic curves over all finite extensions of $k$, implies the section conjecture for all hyperbolic DM-curves over $k$.
\end{theoremB}
When $\U=\Ell_{k}$ is the moduli stack of elliptic curves over $k$, every section $s \in \mathscr{S}_{\pi_1(\Ell_k/k)}$ admits an automorphism of order 2, corresponding to the fact that every elliptic curve has $-1$ as an automorphism. This allows us to \textit{twist} sections $s \in \mathscr{S}_{\pi_1(\Ell_k/k)}$ by an arbitrary quadratic character $\epsilon \colon \Gal_{k} \to \{\pm1\}$, producing a quadratic twist $s\otimes \epsilon\in \mathscr{S}_{\pi_1(\Ell_k/k)}$.

\vspace{\baselineskip}

The $\l$-torsion of the universal elliptic curve over $\Ell_{k}$ yields a monodromy representation 
\[\overline{\rho}_\l \colon \pi_1(\Ell_{k})\longrightarrow \GL_2(\F_\l)\]
and we study $\overline{\rho}_{s,\l}:=\overline{\rho}_\l \circ s$ for (Selmer) sections $s \in  \mathscr{S}_{\pi_1(\Ell_k/k)}$. As one would expect, the representation $\overline{\rho}_{s,\l}$ is compatible with quadratic twists  (see Lemma \ref{lem:quadratictwistcompatibility}), we have 
\[\overline{\rho}_{s\otimes \epsilon,\l}\iso \overline{\rho}_{s,\l} \otimes\epsilon.\]

Assuming the full faithfulness in the section conjecture holds over all finite extensions of the base-field, (quadratic) twisting has the crucial property of preserving rational and cuspidal sections, respectively. Since we work with Selmer sections, it is essential that this full faithfulness holds over number fields and finite extensions of $\Q_p$. Building once again on Bresciani's work, we extend \cite[Theorem~7.2(1)]{brescianiImplicationsGrothendiecksAnabelian2021} from hyperbolic orbicurves to arbitrary hyperbolic DM-curves.
\newtheorem*{theoremC}{Theorem C}
\begin{theoremC}[see Theorem \ref{thm:subpadicisff}]
	Let $k$ be a sub-$p$-adic field, i.e.\ a subfield of a finitely generated extension of $\Q_p$. Then the full faithfulness in the section conjecture holds for all hyperbolic DM-curves over $k$.
\end{theoremC}

\subsection{Outline of the paper} 
Section \ref{sec:dmcurves} introduces DM-curves and is concerned with their étale fundamental groups. In Section \ref{sec:scfordm}, following \cite{brescianiImplicationsGrothendiecksAnabelian2021}, we formulate and study the section conjecture for hyperbolic DM-curves using the étale fundamental gerbe. Section \ref{sec:scforell} treats various properties of sections of the moduli stack of elliptic curves over a general base field. In Section \ref{sec:selmersec}, we study Selmer sections of the moduli stack over $\Q$ via their associated Galois representations and present the proof of the main theorem.

\subsection*{Acknowledgments}
I would like to thank my PhD advisor Jakob Stix for many insightful discussions and comments on earlier drafts. I also thank Magnus Carlson for helpful discussions on anabelian geometry, and Marius Leonhardt for discussions concerning elliptic curves. I further thank Leonie Nienhaus, Ruth Wild and Nils Witt for helpful conversations related to this work. Finally, I would like to thank the anonymous referee for many valuable comments, including the suggestion to add Remark \ref{rem:alternativeapproach}.

The present article is part of the PhD thesis of the author. This project was funded by Deutsche
Forschungsgemeinschaft (DFG, German Research Foundation) through the Collaborative Research
Centre TRR 326 \textit{Geometry and Arithmetic of Uniformized Structures} - Project-ID 444845124.
\subsection*{Notation}
The letter $k$ will denote a field of characteristic zero with algebraic closure $\bar{k}$.  A \textbf{curve} over $k$ is a connected, one-dimensional scheme, smooth and separated over $k$. Finite \'etale morphisms of DM-stacks $ \V\to \U$ are always assumed to be representable and such a morphism is called a \textbf{finite \'etale cover} if it is surjective. For a DM-stack $\U$ and a geometric point $\bar{x} \in \U$, we denote by $\pi_1(\U,\bar{x})$ the \'etale fundamental group, see \cite{noohiFundamentalGroupsAlgebraic2004} for a definition for stacks. When the choice of a base-point is of no consequence, we simply write $\pi_1(\U)=\pi_1(\U,\bar{x})$. For a stack $\U$ over $k$, we write $\overline{\U}$ for the base-change $\U_{\bar{k}}$.  
\section{DM-curves and their fundamental groups}\label{sec:dmcurves}
\subsection{DM-curves} We consider the following generalization of a curve.
\begin{defi}
	 A \textbf{DM-curve} over $k$ is a connected, one-dimensional Deligne-Mumford stack $\U/k$, such that $\U \to \spec(k)$ is smooth and separated. 
\end{defi}
We begin by recalling the following well-known consequence of Chow’s Lemma for algebraic spaces \cite[Theorem 7.4.1]{olssonAlgebraicSpacesStacks2023}.
\begin{lem}\label{lem:algebraicspaceiscurve}
	A one-dimensional algebraic space, which is smooth and separated over a field is a scheme.
\end{lem}
A not necessarily representable morphism of DM-stacks is called \textbf{finite}, if it is proper and quasi-finite.
\begin{lem}
	Let $\U$ be a DM-curve over $k$.
	\be 
	\item The coarse moduli space $j\colon \U \to U$ exists, $j$ is finite and $U/k$ is a curve.
	\item  The DM-curve  $\U/k$ is proper if and only if $U/k$ is proper.
	\ee 
\end{lem}
\bp 
(1) Since $\U/k$ is a separated DM-stack, the coarse moduli space $j\colon \U \to U$ exists and $j$ is proper and quasi-finite by the Keel-Mori Theorem  \cite[Theorem (6.12)]{rydhExistencePropertiesGeometric2013}. By \cite[Lemma 2.2.3]{abramovichCompactifyingSpaceStable2002} \'etale locally on $U$, the morphism $\U \to U$ is of the form 
\[ [V/G] \to V/G,\]
where $V$ is a (not necessarily connected) curve and $G$ a finite constant group. In particular, $U/k$ is separated, one-dimensional and smooth. Hence $U/k$ is a scheme by Lemma \ref{lem:algebraicspaceiscurve} and in particular a curve. \\
(2) If $U/k$ is proper, then so is $\U/k$, because $j$ is proper. If $\U/k$ is proper, then $U/k$ is proper by \cite[Proposition 10.1.6 ($v$)]{olssonAlgebraicSpacesStacks2023}.
\ep 
Recall that a DM-stack is called \textbf{uniformizable} \cite[Definition 6.1]{noohiFundamentalGroupsAlgebraic2004}, if it admits a finite \'etale cover by an algebraic space. For DM-curves, being uniformizable implies that there is a finite \'etale cover by a curve.
\begin{defilem}\label{lem:curvecover}
	Let $\U$ be a DM-curve over $k$. A \textbf{curve-cover} $V\to\U$ is a finite \'etale cover such that $V/k$ is a curve.	The following are equivalent.
	\begin{enumerate}[label=(\alph*)]
		\item We have $\U \iso [V/G]$ for some curve $V/k$ and a finite group $G$ acting on $V/k$.
		\item There is a Galois curve-cover $V\to \U$.
		\item There is a curve-cover $V\to \U$.
		\item The DM-curve $\U$ is uniformizable.
	\end{enumerate}
	In particular, every uniformizable DM-curve admits a Galois curve-cover.
\end{defilem}
\bp 
It is enough to show that (d) implies (a). Assume that $\U$ is uniformizable. Let $f\colon V\to \U$ be a finite \'etale cover by an algebraic space. We may assume that $V$ is connected by replacing it with a connected component. Since $f$ is finite \'etale, $V$ is smooth and separated over $k$. By Lemma \ref{lem:algebraicspaceiscurve} we conclude that $V$ is a curve over $k$. After replacing $V\to\U$ with its Galois closure, we obtain a $G$-Galois cover $V\to\U$ for some finite group $G$. Therefore $\U \iso [V/G]$.
\ep 
The \textbf{Euler-characteristic} of a uniformizable DM-curve $\U$ is defined as 
\[\chi(\U):=\frac{\chi(V)}{\deg(V/\U)},\]
where $V\to \U$ is some curve-cover and $\chi(V)$ is the  $\l$-adic \'etale Euler-characteristic with compact support. Since $k$ is always assumed to be of characteristic zero, the Euler-characteristic $\chi(\U)$ is independent of the choice of the curve-cover $V\to \U$ and multiplicative with respect to connected finite \'etale covers. We call $\U$ \textbf{hyperbolic}, if $\chi(\U)<0$.
\subsection{Fundamental groups}
Recall that given a geometric point $\bar{x}$ of a DM-stack $\U$ and a finite \'etale morphism $\V\to \U$, the group $\Aut(\bar{x})$ acts on the underlying set of the geometric fiber $\V_{\bar{x}}$ (see \cite[page 78]{noohiFundamentalGroupsAlgebraic2004}). This action is functorial in $\V$ and defines a canonical homomorphism 
\[\omega_{\bar{x}} \colon \Aut(\bar{x})\longrightarrow \pi_1(\U,\bar{x}).\]

\begin{thm}\label{thm:noohiuni}
	Let $\U$ be a DM-curve over $k$ with coarse moduli space $j\colon \U \to U$.
	\be 
	\item The induced morphism $\pi_1(j)\colon \pi_1(\U)\twoheadrightarrow \pi_1(U)$ is surjective and the  kernel is the normal subgroup generated by the images of $\omega_{\bar{x}}$  for all $\bar{x} \in \U$.
	\item  The DM-curve $\U$ is uniformizable if and only if the homomorphism 
	\[\omega_{\bar{x}} \colon \Aut(\bar{x})\longrightarrow \pi_1(\U,\bar{x})\]
	is injective for all geometric points $\bar{x} \in \U$.
	\ee  
\end{thm}
\bp 
(1) and (2) are special cases of \cite[Theorem 7.1]{noohiFundamentalGroupsAlgebraic2004} and \cite[Theorem 6.2]{noohiFundamentalGroupsAlgebraic2004}.
\ep 

Associated to a geometrically connected DM-curve $\U/k$, we have the fibration sequence 
\[1 \to \pi_1(\overline{\U})\to \pi_1(\U) \to \Gal_{k} \to 1\]
with $\overline{\U}=\U_{\bar{k}}$. Since any geometric point $\bar{x} \in \mathcal{U}$ lifts to a geometric point of $\overline{\mathcal{U}}$, the functoriality of $\omega_{\bar{x}}$ (see \cite[Lemma 5.1]{noohiFundamentalGroupsAlgebraic2004}) implies that the image of $\omega_{\bar{x}}$ is contained in the subgroup $\pi_1(\overline{\mathcal{U}}) \subset \pi_1(\mathcal{U})$.

\begin{prop}\label{prop:torsionfreeiffcurve}
		Let $\U$ be a uniformizable DM-curve over $k$.
\be 
\item If $k$ is algebraically closed, then $\pi_1(\U)$ is virtually torsion-free and (topologically) finitely presented. 
\item Assume that $\U/k$ is geometrically connected. Then $\U/k$ is a curve if and only if $\pi_1(\overline{\U})$ is torsion-free.
\ee 
\end{prop}
\bp 
(1) Assume that $k$ is algebraically closed. Let $V\to \U$ be a $G$-Galois curve-cover for some finite group $G$. Then we have an exact sequence
\[1 \to \pi_1(V)\to \pi_1(\U)\to G \to 1.\] 
Therefore $\pi_1(\U)$ is an extension of finitely presented groups, hence finitely presented \cite[Claim 2.7 (i)]{esnaultFinitePresentationTame2022}. Since $\pi_1(V)$ is torsion-free, $\pi_1(\U)$ is virtually torsion-free.
(2) When $\U/k$ is a curve, then it's well-known that $\pi_1(\overline{\U})$ is torsion-free. Assume that $\pi_1(\overline{\U})$ is torsion-free. By Lemma \ref{lem:algebraicspaceiscurve} it is enough to show that $\U$ is an algebraic space. By \cite[Theorem 2.2.5 (1)]{conradARITHMETICMODULIGENERALIZED2007} $\U$ is an algebraic space, if the automorphism functors of all geometric points are trivial. Let $\bar{x}\in \overline{\U}$ be a geometric point. By Theorem \ref{thm:noohiuni}, the finite group $\Aut(\bar{x})$ injects into $\pi_1(\overline{\U})$. Since $\pi_1(\overline{\U})$ is torsion-free, this implies  $\Aut(\bar{x})=1$.
\ep 

\subsection{Inertia and decomposition groups}
We extend the results of \cite[Section 3]{stixCuspidalSectionsAlgebraic2012} to the setting of uniformizable DM-curves.

\vspace{\baselineskip}
Let $\U/k$ be a geometrically connected, uniformizable DM-curve with coarse moduli space $j\colon \U \to U$, smooth completion $U\subset X$ and boundary $Y=X\setminus U$. We call $\mathrm{Cusps}(\U):=Y(\bar{k})$ the set of \textbf{cusps} of $\U$. Given a cusp $y\in \mathrm{Cusps}(\U)$, we denote by $X_y:=\spec(\o_{X,y}^h)$ the henselization of $X$ at $y$. Furthermore we set
\[U_y:= U\times_X X_y \quad \text{and} \quad \U_y:=\U\times_U U_y. \]
Since $(U_y)_{\bar{k}}=\overline{U}_{\! y}$, we have $(\U_y)_{\bar{k}}=\overline{\U}_{\! y}$ and the fibration sequence 
\[1 \to \pi_1( \overline{\U}_{\! y})\to\pi_1(\U_y)\to \Gal_{k(y)} \to 1,\]
is exact, where $k(y)$ is the residue field of $y$.
\begin{lem}\label{lem:structurehenselization}
	Let $y \in \mathrm{Cusps}(\U)$ be a cusp and $\bar{\eta} \in \U_y$ a geometric point. Then  $j\colon \U_y\to U_y$ is the coarse moduli space and the sequence
	\[1 \to \Aut(\bar{\eta}) \xrightarrow{\omega_{\bar{\eta}}} \pi_1(\U_y,\bar{\eta}) \xrightarrow{j}\pi_1(U_y,j(\bar{\eta}))\to 1\] 
	is exact.
\end{lem}
\bp 
The morphism $U_y\to U$ is flat and the formation of coarse moduli spaces are stable under flat base-change \cite[Theorem 11.1.2 (iii)]{olssonAlgebraicSpacesStacks2023}. Therefore the pullback $\U_y\to U_y$ is the coarse moduli space. In particular $\U_y\to U_y$ is a gerbe because $U_y$ has only one point. The sequence being right-exact follows from \cite[Proposition 9.5]{noohiFundamentalGroupsAlgebraic2004}. Left-exactness follows from $\U$ being uniformizable.
\ep 

\begin{rem}\label{rem:completions}
	One may (and we will) work equivalently with completions instead of henselizations. For a cusp $y\in \mathrm{Cusps}(\U)$, let $X^{\wedge}_y:=\spec(\o_{X,y}^{\wedge})$ be the spectrum of the completed local ring at $y$ and define
	\[ U^{\wedge}_y:=U\times_X X^{\wedge}_y \quad \text{and} \quad \U^{\wedge}_y:=\U\times_U U^{\wedge}_y.\] Since $\U_y^{\wedge}\to \U_y$ induces an isomorphism $\pi_1(\U^{\wedge}_y)\lisim \pi_1(\U_y)$ over $\Gal_{k(y)}$, the analogous sequences
	\begin{itemize}
		\item $1\to \pi_1( \overline{\U}^{\wedge}_{\! y})\to \pi_1(\U_y^{\wedge})\to \Gal_{k(y)} \to 1$  and
		\item $1 \to \Aut(\bar{\eta})\to \pi_1(\U_y^{\wedge},\bar{\eta})\to \pi_1(U^{\wedge}_y,j(\bar{\eta}))\to 1$
		
	\end{itemize}
	are exact.
\end{rem}

We fix a universal pro-finite \'etale cover $\widetilde{\U}=\lim\limits_{i} \U_i \to \U$, a lift $\widetilde{\U}\to \overline{\U}=\U_{\bar{k}}$ and an isomorphism of the extension
\[1 \to \pi_1(\overline{\U})\to \pi_1(\U)\to \Gal_{k}\to 1\] with the opposite groups of 
\[ 1 \to \Aut(\widetilde{\U}/\overline{\U}) \to \Aut(\widetilde{\U}/\U) \xrightarrow{\mathrm{pr}} \Aut(\spec(\bar{k})/\spec(k)) \to 1.\]
Since $\U$ is uniformizable, we may assume that each finite \'etale cover $\U_i\to \U$ is a curve-cover. Let $\widetilde{\U} \subset \widetilde{\X}$ be the   ``universal ramified cover'' defined by the pro-system induced by the smooth completions of the curves $\U_i$.
Then $\pi_1(\U)$ acts on $\widetilde{\X}$ and in particular on $\widetilde{\Y}:=(\widetilde{\X}\setminus\widetilde{\U})^{\mathrm{red}}$. The (pro-finite) set of \textbf{prolongations of cusps} to $\widetilde{\U}$ is \[\widetilde{\mathrm{Cusps}}(\U):= \Hom_{\bar{k}}(\spec(\bar{k}), \widetilde{\Y}),\]  which carries a natural $\pi_1(\U)$-action given by 
\[ \gamma . \tilde{y} :=  \gamma^{-1} \circ \tilde{y} \circ \pr(\gamma).\]
We have a natural map $\widetilde{\X}\to X$, which restricts to a map $\widetilde{\Y} \to Y$. The induced map 
\[  \widetilde{\mathrm{Cusps}}(\U) \longrightarrow  \Hom_{\bar{k}}(\spec(\bar{k}),Y)=\mathrm{Cusps}(\U)  \]
is $\pi_1(\U)$-equivariant. The \textbf{inertia group} $I_{\tilde{y}/y}$ (resp.\ the \textbf{decomposition group} $D_{\tilde{y}/y}$) is the stabilizer of $\tilde{y}$ under the action of $\pi_1(\overline{\U})$ (resp.\ $\pi_1(\U)$).

\begin{prop}\label{prop:henselizationiso}
Let $\U/k$ be a geometrically connected, uniformizable DM-curve,  $y \in \mathrm{Cusps}(\U)$ a cusp and $\bar{\eta}$ a geometric generic point of $\U$.
\be 
\item A lift $\bar{\eta}_y$ of $\bar{\eta}$ to $\U_y$ determines a prolongation $\tilde{y}$ of $y$ and the image of the homomorphism $\pi_1(\U_y,\bar{\eta}_y)\to \pi_1(\U,\bar{\eta})$ is $D_{\tilde{y}/y}$.
\item If $\chi(\U)\leq 0$, then $\U_y\to \U$ induces an isomorphism of exact sequences  \[\begin{tikzcd}
	1 & {\pi_1(\overline{\U}_{\! y})} & {\pi_1(\U_y)} & {\Gal_{k(y)}} & 1 \\
	1 & {I_{\tilde{y}/y}} & {D_{\tilde{y}/y}} & {\Gal_{k(y)}} & 1.
	\arrow[from=1-1, to=1-2]
	\arrow[from=1-2, to=1-3]
	\arrow["\sim", from=1-2, to=2-2]
	\arrow[from=1-3, to=1-4]
	\arrow["\sim", from=1-3, to=2-3]
	\arrow[from=1-4, to=1-5]
	\arrow[no head, from=1-4, to=2-4]
	\arrow[shift left, no head, from=1-4, to=2-4]
	\arrow[from=2-1, to=2-2]
	\arrow[from=2-2, to=2-3]
	\arrow[from=2-3, to=2-4]
	\arrow[from=2-4, to=2-5]
\end{tikzcd}\] 
\ee 
\end{prop}
\bp 
(1) For any curve-cover $V\to \U$, the lift $\bar{\eta}_y$ of $\bar{\eta}$ determines a geometric point in $V\times_\U \U_y$, which determines a cusp $y'$ of $V$ above $y$. Applying this to $\widetilde{\U}= \lim_i \U_i \to \U$ yields a prolongation $\tilde{y}$ of $y$. For any cusp $y'$ of $\U_i\to \U$ in $\tilde{y}$, a deck transformation in the image of
\[\pi_1(\U_y,\bar{\eta}_y)\to \pi_1(\U,\bar{\eta}) \]
must map the punctured neighborhood $V_{y'}$ to itself and therefore fix the cusp $y'$. This implies that 
\[\im(\pi_1(\U_y,\bar{\eta}_y)\to \pi_1(\U,\bar{\eta}) )\subset D_{\tilde{y}/y}.\] 
Conversely an element in $D_{\tilde{y}/y}$ induces an automorphism of the pro-finite \'etale cover \[\lim_{i}\U_{i,y'}\to \U_y,\] thus lifts to an element of $\pi_1(\U_y,\bar{\eta}_y)$.

(2) Due to (1) it remains to prove injectivity. We may assume that $k=\bar{k}$. By Lemma \ref{lem:structurehenselization}, the sequence
\[1 \to \Aut(\bar{\eta})\to  \pi_1(\U_y,\bar{\eta}) \to  \pi_1(U_y,j(\bar{\eta})) \to  1 \]
is exact.  Let $\gamma \in \ker(\pi_1(\U_y)\to \pi_1(\U))$. If $\gamma$ is of finite order, then $\gamma$ is contained in $\Aut(\bar{\eta})$ as $ \pi_1(U_y)$ is torsion-free. Since $\U$ is uniformizable, $\Aut(\bar{\eta}) \to\pi_1(\U)$ is injective, hence $\gamma=1$. If $\gamma$ is not of finite order, we may replace $\gamma$ by a suitable power  $\gamma^n$  and assume that $\gamma$ is contained in $\pi_1(V_{y'})\subset \pi_1(\U_y)$ for some cusp $y'$ of a curve-cover $V\to \U$, necessarily satisfying $\chi(V)\leq 0$. By \cite[Lemma 14.]{stixCuspidalSectionsAlgebraic2012} the homomorphism  $\pi_1(V_{y'})\to \pi_1(V)\subset \pi_1(\U)$ is injective. Hence $\gamma=1$, which is a contradiction.
\ep

\section{The section conjecture for hyperbolic DM-curves}\label{sec:scfordm}
In this section we formulate and study the section conjecture for hyperbolic DM-curves. We follow \cite{brescianiImplicationsGrothendiecksAnabelian2021} and use the \'{e}tale fundamental gerbe $\Pi_{\mathcal{U}/k}$, which has the advantage of being base-point free and thus enjoys simpler functorial behavior. For an introduction to the \'etale fundamental gerbe see \cite{borneNoriFundamentalGerbe2014} and \cite[Appendix]{brescianiImplicationsGrothendiecksAnabelian2021}. 

\vspace{\baselineskip}

The \'{e}tale fundamental gerbe $\Pi_{\mathcal{U}/k}$ exists if and only if $\U/k$ is geometrically connected \cite[Theorem A.17]{brescianiImplicationsGrothendiecksAnabelian2021}. Accordingly, we henceforth assume that all hyperbolic DM-curves over $k$ are geometrically connected.

\subsection{The section conjecture for hyperbolic DM-curves}
Let $\U/k$ be a hyperbolic DM-curve with coarse moduli space $j\colon \U \to U$, smooth completion $U\subset X$ and boundary $Y=X\setminus U$. The \textbf{rational Kummer map} is the canonical morphism $\U\to \Pi_{\U/k}$ evaluated at $k$-points 
\[\kappa_{\U/k}^\rat\colon \U(k)\longrightarrow \Pi_{\U/k}(k).\]
For a $k$-rational cusp $y\in Y(k)$, the map $\U_y\to \U$ induces a morphism $\Pi_{\U_y/k}\to \Pi_{\U/k}$. Taking the disjoint union over all such $y$ and evaluating at $k$-points, yields the \textbf{cuspidal Kummer map}
\[\kappa^{\cusp}_{\U/k}\colon \bigsqcup_{y \in Y(k)}\Pi_{\U_y/k}(k)\longrightarrow  \Pi_{\U/k}(k).\]
Altogether we obtain the \textbf{non-abelian Kummer map}
\[\kappa_{\U/k}\colon \U(k)\sqcup \bigsqcup_{y \in Y(k)}\Pi_{\U_y/k}(k)\longrightarrow \Pi_{\U/k}(k).\]
Sections \( s \in \Pi_{\mathcal{U}/k}(k) \) in the essential image of $\kappa_{\U/k}^\rat$ (resp.\ $\kappa_{\U/k}^\cusp$) are called  \textbf{rational} (resp.\ \textbf{cuspidal}).

\vspace{\baselineskip}
 For a short exact sequence of profinite groups
 \[
 1 \to \bar{\pi} \to \pi \to \Gamma \to 1,
 \]
 let $\mathrm{Sec}(\pi/\Gamma)$ denote the groupoid whose objects are sections
 $s\colon \Gamma \to \pi$ and whose morphisms $s\to t$ are the elements
 $\gamma\in\bar{\pi}$ satisfying $t=\gamma s\gamma^{-1}$. In particular, 
 \[
 \Aut_{\mathrm{Sec}(\pi/\Gamma)}(s)= \{\gamma \in \bar{\pi} \mid \gamma s \gamma^{-1}=s\}=:Z(s),
 \]
 where $Z(s)=C_{\bar{\pi}}(s(\Gamma))$ is the centralizer of the image of $s$. By \cite[Proposition 9.3]{borneNoriFundamentalGerbe2014}, we have equivalences of groupoids
\[
\Pi_{\U/k}(k) \simeq \mathrm{Sec}(\pi_1(\U,\bar{x})/\Gal_k) \quad  \text{and} \quad 
\Pi_{\U_y/k}(k) \simeq \mathrm{Sec}(\pi_1(\U_y,\bar{\eta})/\Gal_k)
\]
for any choice of geometric points $\bar{x} \in \overline{\U}$ and $\bar{\eta} \in \overline{\U}_{\! y}$. In particular the automorphisms of an object $s \in \Pi_{\U/k}(k)$ are identified with the centralizer $Z(s)$ of the corresponding section. See \cite[Lemma 2.6.]{brescianiImplicationsGrothendiecksAnabelian2021}, where this identification is explained in detail. We denote by $\mathscr{S}_{\pi_1(\U/k)}$ (resp.\ $\mathscr{S}_{\pi_1(\U_y/k)}$) the respective sets of conjugacy classes of sections.
In particular, when $\U = U$ is a curve, the non-abelian Kummer map $\kappa_{\U/k}$ provides a categorical lift of the classical non-abelian Kummer map.
\begin{defi}
	Let $\U/k$ be a hyperbolic DM-curve. We say that:
	\be 
	\item  The  \textbf{(stacky) section conjecture} holds for $\U/k$, if $\kappa_{\U/k}$ is an equivalence.
	\item The \textbf{surjectivity in the section conjecture} holds for $\U/k$, if $\kappa_{\U/k}$ is essentially surjective.
	\item The  \textbf{injectivity in the section conjecture}  holds for $\U/k$, if $\pi_0(\kappa_{\U/k})$ is injective.
		\item The  \textbf{full faithfulness in the section conjecture}  holds for $\U/k$, if $\kappa_{\U/k}$ is fully faithful.
	\ee 
\end{defi}
\begin{defi}
	Let $\U/k$ be a hyperbolic DM-curve. A \textbf{curve-neighborhood} of a section $s \in \Pi_{\U/k}(k)$ consists of a geometrically connected curve-cover $f\colon V\to \U$ together with a section $\tilde{s} \in \Pi_{V/k}(k)$, such that $f(\tilde{s})\iso s$.
\end{defi}
We will use the following lemma repeatedly to reduce to the case of curves.
\begin{lem}\label{lem:curveneighborhood}
 Let $\U/k$ be a hyperbolic DM-curve and $s \in \Pi_{\U/k}(k)$ a section. Then $s$ has a curve-neighborhood $V\to \U$, such that $\pi_1(\overline{V}) \subset \pi_1(\overline{\U})$ is a characteristic subgroup.
\end{lem}
\bp 
For $m\geq 1$, consider the descending sequence of characteristic subgroups 
\[N_m:=\bigcap_{\phi \colon \pi_1(\overline{\U}) \to G}\ker{\phi}\]
where $G$ ranges over all finite groups with $\# G\leq m$. By Proposition \ref{prop:torsionfreeiffcurve} (1) the group $\pi_1(\overline{\U})$ is finitely generated and hence $N_m \subset \pi_1(\overline{\U})$ has finite index. Moreover, the family of subgroups $(N_m)_m$ is cofinal in the set of all finite-index subgroups of $\pi_1(\overline{\U})$. Since $\U$ is uniformizable, there exists a sufficiently large $m_0$, such that $\overline{N}:=N_{m_0}$ is the fundamental group of a curve-cover $\overline{V}\to \overline{\U}$. Since $\overline{N}$ is characteristic, the subgroup
\[ N:=\overline{N}s(\Gal_{k}) \subset \pi_1(\U)\]
satisfies $N\cap \pi_1(\overline{\U})=\overline{N}$. Thus, the corresponding connected finite \'etale cover $V\to \U$  satisfies $V_{\bar{k}} \iso \overline{V}$ and is a curve-neighborhood of $s$.
\ep 
Using Lemma \ref{lem:curveneighborhood} we obtain the following.
\begin{lem}\label{lem:surjectivityimpliesessential}
	Assume that the surjectivity in the section conjecture holds for all hyperbolic curves over $k$. Then the surjectivity in the section conjecture holds for all hyperbolic DM-curves over $k$.
\end{lem}
\begin{lem}\label{lem:goingupanddown}
	Let $f\colon \V\to \U$ be a finite \'etale cover of (geometrically connected) hyperbolic DM-curves over $k$. Let $t \in \Pi_{\V/k}(k)$ be a section and $s:=f(t) \in \Pi_{\U/k}(k)$. Then $t$ is cuspidal (resp.) rational if and only if $s$ is cuspidal (resp.) rational.
\end{lem}
\bp 
If $t$ is cuspidal (resp.) rational, then $s$ is cuspidal (resp.) rational by the functoriality of the non-abelian Kummer map. Assume that $s$ is cuspidal. Then $s(\Gal_k) \subset D_{\tilde{y}/y}$ for some prolongation $\tilde{y}$ of a $k$-rational cusp $y$ of $\U$. Let $y'$ be the corresponding cusp of $\V$. Hence 
\[t(\Gal_k) \subset D_{\tilde{y}/y} \cap \pi_1(\overline{\V}) = D_{\tilde{y}/y'}. \]
Hence $y'$ is $k$-rational and $t$ lies in the image of 
\[ \Pi_{\V_{y'}/k}(k) \longrightarrow \Pi_{\V/k}(k).\]
Therefore $t$ is cuspidal. Assume that $s$ is rational, then $s$ is the image of $\kappa_{\U/k}^\rat$. By \cite[Proposition A.25]{brescianiImplicationsGrothendiecksAnabelian2021}, we have
 \[\V(k) \iso \U(k)\times _{\Pi_{\U/k}(k)} \Pi_{\V/k}(k),\]
 which implies that $t$ lies in the image of $\kappa_{\V/k}^\rat$. Hence $t$ is rational.
\ep 
We have the following generalization of \cite[Theorem 14]{stixCuspidalSectionsAlgebraic2012} and \cite[Theorem 7.2 (2)]{brescianiImplicationsGrothendiecksAnabelian2021}.
\begin{thm}\label{thm:abstractff}
	Let $k$ be a field of characteristic zero.
	\be
		\item The injectivity in the section conjecture for all proper hyperbolic curves over all finite extensions of $k$, implies the full faithfulness in the section conjecture for all hyperbolic DM-curves over $k$.
		\item The section conjecture for all proper hyperbolic curves over all finite extensions of $k$, implies the stacky section conjecture for all hyperbolic DM-curves over $k$.
	\ee
\end{thm}
\bp 
(1) Step 1: We claim that  $\kappa_{V/k}^\rat$ is fully faithful for all hyperbolic curves $V/k$.

 Let $x,y \in V(k)$, such that $s_x=s_y$ in $\mathscr{S}_{\pi_1(V/k)}$. After passing to a neighborhood of $s_x=s_y$, we may assume that the smooth completion $i\colon V \hookrightarrow X$ is hyperbolic. Then \[s_{i(x)}=i(s_x)=i(s_y)=s_{i(y)}\] in $\mathscr{S}_{\pi_1(X/k)}$. The injectivity in the section conjecture for $X/k$ implies $x=y$. Hence $\pi_0(\kappa_{V/k}^\rat)$ is injective. It remains to prove that $Z(s_x)=1$. We adapt the proof in \cite[Proposition 3.1]{brescianiImplicationsGrothendiecksAnabelian2021}. Let $\phi\colon \pi_1(\overline{V})\twoheadrightarrow \overline{G}$ be a characteristic finite quotient with corresponding cover $\overline{W}\to \overline{V}$. For $\gamma \in Z(s_x)$ we set $g:=\phi(\gamma)$. Using $x \in V(k)$, we may descend $\overline{W}\to \overline{V}$ to a cover $W\to V$ lifting $x$, which is a torsor under a $k$-form $G$ of $\overline{G}$. Let $w \in W(k)$ be a point over $x$. Since $g$ is $\Gal_k$-invariant, we have $g \in G(k)$ and $w':=g.w \in W(k)$.  We view $g \in \overline{G}=\Aut(\overline{W}/\overline{V})$ as a deck transformation. Then $g$ and conjugation  by $\gamma$ induce the same outer automorphism of $\pi_1(\overline{W})$. Hence, up to conjugation by an element in $\pi_1(\overline{W})$, we have
\[s_w = \gamma s_w \gamma^{-1}= g(s_{w})=s_{g(w)}=s_{w'}.\]
By assumption this implies $w=w'$. Since $\overline{G}$ acts freely on $\overline{W}$, we conclude $g=1$. This proves $Z(s_x)=1$ and hence the full faithfulness of $\kappa_{V/k}^\rat$.

Step 2: We claim that  $\kappa_{\U/k}^\rat$ is fully faithful for all hyperbolic DM-curves $\U/k$.

 We may assume $\U(k)\neq \emptyset$, since otherwise there is nothing to prove. Then there is a geometrically connected curve-cover $V\to \U$. By Step 1 the functor $\kappa_{V_{k'}/k'}^\rat$ is fully faithful for all finite extensions $k'/k$. Applying \cite[Proposition 2.5]{brescianiImplicationsGrothendiecksAnabelian2021} to $V\to \U$ proves that $\kappa_{\U/k}^\rat$ is fully faithful.

Step 3: We claim that  $\kappa_{\U/k}^\cusp$ is fully faithful for all hyperbolic DM-curves $\U/k$. 

Let $y,y'$ be $k$-rational cusps of $\U$ and $s \in \mathscr{S}_{\pi_1(\U_y/k)}, t \in \mathscr{S}_{\pi_1(\U_{y'}/k)}$ sections, such that $s= t$ in $\mathscr{S}_{\pi_1(\U/k)}$. By Lemma \ref{lem:goingupanddown}, we may pass to a curve-neighborhood and assume that $\U=V$ is a hyperbolic curve with hyperbolic smooth completion $i \colon V\hookrightarrow X$. Then 
\[s_y = i(s)= i(t) =s_{y'}\]
in $ \mathscr{S}_{\pi_1(X/k)}$ and the injectivity in the section conjecture for $X/k$ implies $y=y'$. Thus, it suffices to prove that the functor
\[\Pi_{\U_y/k}(k)\to \Pi_{\U/k}(k) \]
is fully faithful. First, notice that it is faithful. This follows from the injectivity of the map $\pi_1(\U_y)\to \pi_1(\U)$  by Proposition \ref{prop:henselizationiso} (2). Let $s,t  \in \Pi_{\U/k}(k)$ and $\gamma \in \pi_1(\overline{\U})$, such that $ s =\gamma t \gamma^{-1}$. Passing over all curve-neighborhoods, the sections $s$ and $t$ determine unique prolongations $\tilde{y},\tilde{y}'$ of $y$, satisfying
\[ s(\Gal_{k})\subset D_{\tilde{y}/y} \text{ and } t(\Gal_{k})\subset D_{\tilde{y}'/y}.\]
Due to the uniqueness, we find 
\[ D_{\tilde{y}/y}= \gamma D_{\tilde{y}'/y} \gamma^{-1}= D_{\gamma.\tilde{y}'/\gamma.y}, \]
which implies $\gamma.y= y$. Applying this argument to all neighborhoods, we find $\gamma.\tilde{y}=\tilde{y}$ and hence $\gamma \in I_{\tilde{y}/y}= \pi_1(\U_y)$.

Step 4: Let $\U/k$ be a hyperbolic DM-curve over $k$. We claim that $\kappa_{\U/k}$ is fully faithful.

 By Step 2 and 3, it remains to prove that cuspidal and rational sections are disjoint. Let $s \in \mathscr{S}_{\pi_1(\U/k)}$ be a section which is cuspidal and rational at the same time. After passing to a curve-neighborhood of $s$, we may assume that $\U=V$ is a hyperbolic curve with hyperbolic smooth completion $i \colon V \hookrightarrow X$. But then $i(s) \in \mathscr{S}_{\pi_1(X/k)}$ is rational and associated to a point in $V(k)$ and in $(X\setminus V)(k)$ at the same time. This contradicts the injectivity in the section conjecture for $X/k$.
 
(2) Because of (1), it suffices to prove that the surjectivity in the section conjecture holds for hyperbolic DM curves over $k$. Let $\U/k$ be a hyperbolic DM-curve and $s \in \mathscr{S}_{\pi_1(U/k)}$.  Using Lemma \ref{lem:surjectivityimpliesessential} we find a curve-cover $V\to \U$ with lift $\tilde{s} \in\mathscr{S}_{\pi_1(V/k)}$. If $\tilde{s}$ is cuspidal or rational, then so is $s$, thus we may replace $\U$ with $V$ and assume that $\U=V$ is a hyperbolic curve.  After further replacing $V$ with a neighborhood of $s$, we may assume that the smooth completion $i\colon V \hookrightarrow X$ is hyperbolic. By assumption $i(s) \in  \mathscr{S}_{\pi_1(X/k)}$ is associated to a unique rational point $x \in X(k)$. Applying this argument to all neighborhoods $V'\to V$ of $s$ with smooth completions $V'\subset X'$, we obtain a unique point
\[\tilde{x} \in \lim_{V'} X'(k).\]
If $x$ is a cusp of $V$, then $\tilde{x}$ is a prolongation of $x$ satisfying $s(\Gal_k)\subset D_{\tilde{x}/x}$. Hence $s$ is cuspidal. If $x \in V(k)$ is a rational point, the same applies to all neighborhoods and we have 
\[\tilde{x} \in \lim_{V'} V'(k).\]
Thus $s=s_x$ is by \cite[Lemma 52.]{stixRationalPointsArithmetic2012}, which implies that $s$ is rational.
\ep 
\begin{rem}\label{rem:alternativeapproach}
	We thank the referee for suggesting that we include the following alternative formulation of the section conjecture for DM-curves. Let $\U/k$ be a geometrically connected, uniformizable DM-curve and $\widetilde{\U}= \lim_i \U_i \to \U$ a universal pro-finite \'etale cover. Since $\U$ is uniformizable, we may assume that all $\U_i$ are curves with smooth completions $X_i$. Let $\widetilde{\X}=\lim_{i}X_i$ be the   ``universal ramified cover''. Because the smooth completion is functorial, the group $\Aut(\widetilde{\U}/\U)=\lim_i \Aut(\U_i/\U)$ acts on $\widetilde{\X}$. We may consider the stack
\[
\mathcal{X}:=\lim_i \;[X_i/\Aut(\U_i/\U)].
\] It is a geometrically connected pro-DM stack over $k$ with an open embedding $\U \subset \mathcal{X}$. The inclusion $\U \subset \mathcal{X}$ induces an isomorphism on fundamental groups, hence an equivalence of \'etale fundamental gerbes 
	\[\Pi_{\U/k} \lisim \Pi_{\mathcal{X}/k}.\]
	One can check that the section conjecture (resp. the full faithfulness in the section conjecture) holds for $\U/k$ if and only if the canonical morphism 
	\[\mathcal{X}(k) \longrightarrow \Pi_{\X/k}(k)\]
	is an equivalence (resp. fully faithful). In the case where $\U=U$ is a curve, this is proven in \cite[Section 8]{brescianiImplicationsGrothendiecksAnabelian2021}, from which the general case may be deduced. 
\end{rem}

		\subsection{Twists and centralizers}
Let $\U/k$ be a geometrically connected DM-curve over a field of characteristic zero  $k$ and $k'/k$ an arbitrary field extension. By \cite[Proposition A.23.]{brescianiImplicationsGrothendiecksAnabelian2021}, we have
\[\Pi_{\U/k}\times_k k' \iso \Pi_{\U_{k'}/k'}.\]
Consequently, the base change of sections is well defined for arbitrary field extensions of $k$.
For a section $s \in \Pi_{\U/k}(k)$, we denote its base-change by $s_{k'} \in  \Pi_{\U_{k'}/k'}(k')$.
\begin{defi}
	Let $\U/k$ be a hyperbolic DM-curve  and $s,t \in \Pi_{\U/k}(k)$ be sections.
	\be 
	\item Let $k'/k$ be a finite extension. We call $t$ \textbf{a twist of $s$ over $k'$}, if $s_{k'} \iso t_{k'}$.
	\item We call $t$ \textbf{a twist of $s$}, if $t$ is a twist of $s$ over $k'$ for some finite extension $k'/k$.
	\ee  
\end{defi}
\begin{lem}\label{lem:beingcuspidaldescends}
	Let $\U/k$ be a hyperbolic DM-curve and $s \in \Pi_{\U/k}(k)$.  Assume that for all finite extensions $k'/k$, the full faithfulness in the section conjecture holds for $\U_{k'}/k'$. 
	\be 
	\item If $s_{k_0}$ is cuspidal (resp.\ rational) for some finite extension $k_0/k$, then $s$ is cuspidal (resp.\ rational). 
	\item Let $t\in \Pi_{\U/k}(k)$ be a twist of $s$. Then $t$ is cuspidal (resp.\ rational) if and only if $s$ is cuspidal (resp.\ rational).
	\ee 
\end{lem}
\bp 
(2) is a direct consequence of (1), hence it suffices to prove (1). Let $j\colon \U \to U$ be the coarse moduli space with smooth completion $U\subset X$ and boundary $Y=X\setminus U$. By assumption the canonical morphisms
\[\U \longrightarrow \Pi_{\U/k} \quad \text{and} \bigsqcup_{y \in Y(\bar{k})}\Pi_{\U_y/k(y)}\longrightarrow \Pi_{\U/k}\]
are fully faithful as morphisms of stacks over the small \'etale site of $\spec(k)$. In particular, the essential images are substacks of  $\Pi_{\U/k}$. Hence membership in the respective essential images can be checked \'etale locally on $\spec(k)$.
\ep 

Let $\U/k$ be a hyperbolic DM-curve and $s \in \Pi_{\mathcal U/k}(k)$ a section. Fix a geometric point $\bar{x} \in \overline{\mathcal U}$, and identify $s$ with the corresponding section of the exact sequence
\[
1 \longrightarrow \pi_1(\overline{\mathcal U}, \bar{x})
\longrightarrow \pi_1(\mathcal U, \bar{x})
\longrightarrow \Gal_k
\longrightarrow 1 .
\]
Then  $\Aut(s) \simeq Z(s)$, where $Z(s)$ denotes the centralizer of $s(\Gal_k)$ in $\pi_1(\overline{\mathcal U}, \bar{x})$. Moreover, for a finite Galois extension $k'/k$, the isomorphism classes of twists of $s$ over $k'$ are classified by the non-abelian cohomology set $\H^1(k'/k,\, Z(s_{k'}))$.
Define
\[
\underline{Z}(s) := \bigcup_{k'} Z(s_{k'}) \subset \pi_1(\overline{\mathcal U}, \bar{x}),
\]
where the union ranges over all finite Galois extensions $k'/k$. Then $\underline{Z}(s)$ is a discrete $\Gal_{k}$-group and the set 
\[\H^1\bigl(k,\, \underline{Z}(s)\bigr)=\colim_{k'/k}\H^1(k'/k,\, Z(s_{k'}))\] classifies arbitrary twists of $s$ over $k$ up to isomorphism.

\vspace*{\baselineskip}

Let $x\in\U(k)$ and assume that the full faithfulness in the section conjecture holds for hyperbolic DM-curves over $k$. Then the rational Kummer map $\kappa_{\U/k}^\rat \colon \U(k) \to \Pi_{\U/k}(k)$ induces an isomorphism
\[\Aut(x)\lisim \Aut(s_x)\iso Z(s_x).\]
Thus, if the section conjecture holds for hyperbolic DM-curves over $k$, centralizers of non-cuspidal sections should be finite. For cuspidal sections we have the following.
\begin{lem}\label{lem:centralizerofcuspidal}
	Let $\U/k$ be a hyperbolic DM-curve and let $s \in \Pi_{\U/k}(k)$ be a cuspidal section associated to a $k$-rational cusp $y$ of $\U/k$. Assume that 
	\begin{itemize}
	\item the full faithfulness in the section conjecture holds for $\U/k$ and
	\item $\widehat{\Z}(1)^{\Gal_k}=0$.
	\end{itemize} 
Let $\bar{\eta} \in \overline{\U}_{\! y}$ be a geometric point and identify $s$ with a section of 
\[1 \to \pi_1(\overline{\U}_{\! y},\bar{\eta}) \to \pi_1(\U_y,\bar{\eta})\to \Gal_k\to 1. \]	
	Then $Z(s)\subset \Aut(\bar{\eta})$. 

\end{lem}
\bp 
 The full faithfulness in the section conjecture implies that $Z(s)\subset \pi_1(\overline{\U}_{\! y},\bar{\eta})$. Let $j\colon \U \to U$ be the coarse moduli space. By Lemma \ref{lem:structurehenselization}, the sequence
 \[1 \to \Aut(\bar{\eta})\xrightarrow{\omega_{\bar{\eta}}} \pi_1( \overline{\U}_{\! y},\bar{\eta})\xrightarrow{j}\pi_1(\overline{U}_{\! y},j(\bar{\eta}))\to 1\]
 is exact. The image of $Z(s)$ in $\pi_1(\overline{U}_{\! y},j(\bar{\eta}))=\widehat{\Z}(1)$ is contained in the $\Gal_{k}$-invariants, which are trivial by assumption. Therefore $Z(s)\subset \Aut(\bar{\eta})$.
\ep 
\begin{defi}
	We say that $k$ is \textbf{section-discrete}, if for any hyperbolic curve $V/k$, the groupoid $\Pi_{V/k}(k)$ is discrete, i.e.\ any section has trivial centralizer. We say that $k$ is \textbf{very section-discrete}, if every finite extension $k'/k$ is section-discrete.
\end{defi}
\begin{prop}\label{prop:anabeliancentralizers}
	Let $\U/k$ be a hyperbolic DM-curve and $s \in \Pi_{\U/k}(k)$ a section.
	\be 
	
	\item If $k$ is section-discrete, then the group $Z(s)$ is finite. 
	\item If $k$ is very section-discrete, then the $\Gal_k$-group $\underline{Z}(s)$ is finite.
	\ee 
\end{prop}

\bp 
Assume that $k$ is section-discrete and let $s \in \Pi_{\U/k}(k)$ be a section. Let $V\to \U$ be a $G$-Galois curve-neighborhood of $s$ with lift $\tilde{s}$ for some finite group $G$. Since $k$ is section-discrete, we have
\[ Z(s) \cap \pi_1(\overline{V})=Z(\tilde{s})=1.\]
We conclude that the composition
\[Z(s)=Z(s)/Z(s)\cap  \pi_1(\overline{V}) \hookrightarrow \pi_1(\overline{\U})/\pi_1(\overline{V})\iso G\]
is injective, which implies that $Z(s)$ is finite. If we further assume that $k$ is very section-discrete, then by the same argument $\underline{Z}(s)$ 
injects into $G$. \ep

\begin{prop}\label{prop:semiabelianimpliessectiondiscrete}
	Let $k$ be a field of characteristic zero. Assume that for every semi-abelian variety $B/k$ with Tate-module $TB=\prod_\l T_\l B$, we have $\H^0(k,TB)=0$. Then $k$ is section-discrete.
\end{prop}
\bp 
The proof is exactly the same as \cite[Proposition 104]{stixRationalPointsArithmetic2012}.
\ep 
\subsection{Sub $p$-adic fields and Selmer sections}
Recall that a field $k$ is called sub-$p$-adic, if $k$ is a subfield of a finitely generated extension of $\Q_p$. 
\begin{lem}\label{lem:abelianvarietyoversubpadic}
	Let $k$ be a sub-$p$-adic field and $B/k$ be a semi-abelian variety with Tate-module $TB=\prod_\l T_\l B$. Then $\H^0(k,TB)=0$.
\end{lem}
\bp 
\cite[Lemma 3.2]{brescianiSectionConjectureFields2025}.
\ep 
\begin{thm}\label{thm:subpadicisff}
	Let $k$ be a sub-$p$-adic field. 
	\be 
	\item The full faithfulness in the section conjecture holds for hyperbolic DM-curves over $k$ and
	\item $k$ is very section-discrete. 
	\ee  
\end{thm}
\bp 
 We first note that being sub-$p$-adic is closed under finite extensions. 
 
(1) By \cite[Lemma 3.5]{brescianiSectionConjectureFields2025} the injectivity in the section conjecture holds for all hyperbolic curves over all finite extensions of $k$. Then Theorem \ref{thm:abstractff} implies that the full faithfulness in the section conjecture holds for hyperbolic DM-curves over $k$.

(2) By Lemma \ref{lem:abelianvarietyoversubpadic}, the assumptions of Proposition \ref{prop:semiabelianimpliessectiondiscrete} apply to $k$.
\ep
\begin{rem}
	Part (2) of Theorem \ref{thm:subpadicisff} was also proven in \cite[Lemma 3.3]{brescianiSectionConjectureFields2025}.
\end{rem}
\begin{defi}
	Let $\U$ be a DM-curve over a number field $k$. A section $s \in \Pi_{\U/k}(k)$ is called \textbf{Selmer}, if for all places $v$ of $k$, the base-changed section $s_{k_v} \in \Pi_{\U_{k_v}/k_v}(k_v)$ is cuspidal or rational. The subset of (conjugacy-classes) of Selmer sections is denoted $\mathscr{S}_{\pi_1(\U/k)}^\mathrm{Sel}\subset \mathscr{S}_{\pi_1(\U/k)}$.
\end{defi}
Let $V/k$ be a hyperbolic curve over a number field $k$ with smooth completion $X$. We briefly recall the \textbf{finite descent obstruction} $X(\A_k)^{\mathrm{f\text{-}cov}}$ \cite[Definition 5.4]{stollFiniteDescentObstructions2007}. It is a subset  of the modified set of adelic points  $X(\A_k)_\bullet$ of $X$,  where for infinite places $v$, the set $X(k_v)$ is replaced by the set of connected components $\pi_0(X(k_v))$. An adelic point $(x_v)_v\in X(\A_k)_\bullet$ lies in $X(\A_k)^{\mathrm{f\text{-}cov}}$ if for every torsor $Y\to X$ under a finite \'etale group scheme $G/k$, there exists a class $\xi \in \H^1(k,G)$ such that $(x_v)_v$ lifts to the twist $Y^\xi\to X$.

 The set of Selmer sections $\mathscr{S}_{\pi_1(V/k)}^\mathrm{Sel}$ is intimately related to the set $X(\A_k)^{\mathrm{f\text{-}cov}}$. The local components of a Selmer section $s \in \mathscr{S}_{\pi_1(V/k)}^\mathrm{Sel}$ determine a unique adelic point  $\underline{x}(s)=(x_v)_v \in X(\A_k)_\bullet$ and from \cite[Theorem 11]{harariDescentObstructionFundamental2012} one concludes that $\underline{x}(s) \in X(\A_k)^{\mathrm{f\text{-}cov}}$. Thus, we obtain a map 
\[ \underline{x} \colon \mathscr{S}_{\pi_1(V/k)}^\mathrm{Sel}\to  X(\A_k)^{\mathrm{f\text{-}cov}}. \]
We have the following well-known consequence of Stoll's results on the finite descent obstruction, see \cite[Corollary 5]{stixBirationalSectionConjecture2015}.
\begin{lem}\label{lem:selmeradelic}
Let $V/k$ be a hyperbolic curve over a number field $k$ with smooth completion $X$. A Selmer section $s \in  \mathscr{S}_{\pi_1(V/k)}^\mathrm{Sel}$ is cuspidal or rational if and only if $\underline{x}(s) \in X(k) \subset X(\A_k)^{\mathrm{f\text{-}cov}}$.
\end{lem}
\bp 
If $s$ is rational, associated with a rational point $x \in V(k)$ (resp.\ cuspidal, associated to a cusp $y \in (X \setminus V)(k)$), then $\underline{x}(s) = x \in X(k)$ (resp.\ $\underline{x}(s) = y \in X(k)$).

Conversely, suppose that $\underline{x}(s)=x \in X(k)$. Let $V' \to V$ be a neighborhood of $s$ with lift $\tilde{s} \in  \mathscr{S}_{\pi_1(V'/k)}$ and smooth completion $V' \subset X'$ of genus $\geq 1$. Since $s$ is Selmer, so is $\tilde{s}$. Let $f\colon X'\to X$ be the induced branched cover. The fiber $Z:=f^{-1}(x)\subset X'$ is zero-dimensional and we have $ \underline{x}(\tilde{s}) \in X'(\A_k)^{\mathrm{f\text{-}cov}} \cap Z(\A_k)_\bullet$. By \cite[Theorem 8.2]{stollFiniteDescentObstructions2007}, we have  
\[X'(\A_k)^{\mathrm{f\text{-}cov}} \cap Z(\A_k)_\bullet=Z(k)\]
and hence  $\underline{x}(\tilde{s}) \in X'(k)$. Using Tamagawa's limit argument \cite[Corollary (2.10)]{tamagawaGrothendieckConjectureAffine1997}, we conclude that $s$ is rational or cuspidal. 
\ep 
\begin{prop}\label{prop:finitedescentmoduli}
	Let $\U/k$ be a hyperbolic DM-curve over a number field $k$. Assume there is a non-constant $k$-morphism $\phi \colon \U \to W$ into a proper variety $W/k$ satisfying
\[W(\A_k)^{\mathrm{f\text{-}cov}}=W(k).\]
Then every Selmer section $s \in\mathscr{S}_{\pi_1(\U/k)}^\mathrm{Sel}$ is cuspidal or rational.
\end{prop}
\bp 
Let $s \in \mathscr{S}_{\pi_1(\U/k)}^\mathrm{Sel}$ be a Selmer section and $f \colon V\to \U$ a curve-neighborhood of $s$ with lift $\tilde{s}$ and hyperbolic smooth completion $V\subset X$. By Lemma \ref{lem:goingupanddown}, $\tilde{s}$ is a Selmer section of $V/k$ with adelic point $\underline{x}(\tilde{s}) \in X(\A_k)^{\mathrm{f\text{-}cov}}$. The composition $\phi \circ f$ extends to a non-constant morphism $X \to W$. By \cite[Proposition 8.5]{stollFiniteDescentObstructions2007}, we have $X(\A_k)^{\mathrm{f\text{-}cov}}=X(k)$. Therefore $\underline{x}(\tilde{s}) \in X(k)$ and Lemma \ref{lem:selmeradelic} implies that $\tilde{s}$ is cuspidal or rational. Hence $s=f(\tilde{s})$ is cuspidal or rational.
\ep 
\section{The section conjecture for the moduli stack of elliptic curves}\label{sec:scforell}
Let $k$ be a field of characteristic zero with algebraic closure $\bar{k}$. We denote by $\Ell_{k}/k$ the moduli stack of elliptic curves over $k$. 
\subsection{The moduli stack of elliptic curves as a DM-curve}
\begin{prop}\label{prop:ellisdmcurve}
	The stack $\Ell_{k}$ is a hyperbolic, geometrically connected DM-curve over $k$ with coarse moduli space given by the $j$-invariant $j\colon \Ell_{k}\to \A^1_k$.

\end{prop}
\bp 
The stack $\Ell_{k}$ being a one-dimensional DM-stack, which is smooth and separated over $k$ with moduli space $j\colon \Ell_{k} \to \A^1_k$ is classical. Because $j$ induces a homeomorphism on associated topological spaces and $\A^1_k$ is connected, we conclude that $\Ell_{k}$ is connected. Since $k$ was arbitrary, $\Ell_{k}$ is geometrically connected. The classifying morphism $\P^1_k\setminus \{0,1,\infty\}\to \Ell_{k}$, corresponding to the Legendre family of elliptic curves is finite \'etale \cite[(4.6.2)]{katzArithmeticModuliElliptic1985}, hence a hyperbolic curve-cover.
\ep
By the invariance of the \'etale fundamental group under algebraically closed base change in characteristic zero and \cite[Theorem 1]{odaEtaleHomotopyType1997}, we have
\[
\pi_1(\Ell_{\bar{k}}) \simeq \widehat{\SL_2(\Z)},
\]
with $\widehat{(-)}$ denoting profinite completion.

\begin{lem}\label{lem:sl2}
	The image of $-I \in \SL_2(\Z)\to \widehat{\SL_2(\Z)}$ is the unique element of order 2. In particular it is central.
\end{lem}
\bp 
As $\SL_2(\Z)\subset \widehat{\SL_2(\Z)}$ is dense, the element $-I$ is central. The matrices $S,R \in \SL_2(\Z)$ given by
\[S:=\begin{pmatrix}
	0 & -1 \\
	1 & 0
\end{pmatrix} \quad \text{and} \quad R:=\begin{pmatrix}
	0 & -1 \\
	1 & 1
\end{pmatrix} \]
satisfy $S^2=R^3= -I$ and induce an amalgamated product decomposition 
\[\Z/4\Z \ast_{\Z/2\Z} \Z/6\Z \lisim \SL_2(\Z),\]
see \cite[§1, 1.5.3]{serreArbresAmalgames1977}.
Hence, we have a profinite amalgamated product decomposition 
\[ \Z/4\Z\; \widehat{\ast}_{\Z/2\Z}\; \Z/6\Z \lisim \widehat{\SL_2(\Z)}.\]
 By \cite[Theorem 7.1.2]{ribesProfiniteGraphsGroups2017} any finite subgroup $C\subset \widehat{\SL_2(\Z)}$ is conjugate to a subgroup of $\langle R\rangle $ or $\langle S\rangle$. Therefore any order 2 element is conjugate to $-I$ and hence equal to $-I$.
\ep 
The corresponding element  $-I \in \pi_1(\Ell_{\bar{k}})$ is characteristic, hence lies in the center of $\pi_1(\Ell_k)$. Thus, for any section $s \in \Pi_{\Ell_k/k}(k)$, we have $ \{\pm I\}\subset \underline{Z}(s)$. In particular, we obtain a map
\[ \Hom(\Gal_k,\{\pm 1\})=\H^1(k,\{\pm I\})\to  \H^1(k,\underline{Z}(s)),\]
allowing us to twist a section $s \in \Pi_{\Ell_k/k}(k)$ by a quadratic character 
\[\epsilon\colon \Gal_k\to \{\pm 1\}.\]
The resulting section is called the \textbf{quadratic twist of $s$ by~$\epsilon$} and is denoted $s \otimes \epsilon $.
\begin{rem}
	The natural homomorphism $\widehat{\SL_2(\Z)}\twoheadrightarrow \SL_2(\widehat{\Z})$ is surjective, but not an isomorphism. By a result of Mel'nikov \cite{melnikovCongruenceKernelGroup1976}, the kernel is a free profinite
	group on a countable infinite number of generators.
\end{rem}

\subsection{Cuspidal sections and the Tate curve}
The coarse moduli space of $\Ell_k$ is $\A_k^1$, hence $\Ell_{k}$ has exactly one $k$-rational cusp $\infty \in \P^1(k)$. For a geometric point $\bar{\eta} \in \Ell_{k,\infty}$, the sequence 
\[1 \to \Aut(\bar{\eta})\xrightarrow{\omega_{\bar{\eta}}}\pi_1(\Ell_{k,\infty},\bar{\eta})\xrightarrow{j} \pi_1(\A^1_{k,\infty},j(\bar{\eta}))\to 1\]
is exact by Lemma \ref{lem:structurehenselization}.
Since $\bar{\eta}$ lies over the generic point of $\Ell_k$, we have $\Aut(\bar{\eta}) \iso \Z/2\Z$. In particular, the uniqueness part of Lemma \ref{lem:sl2} implies
\[\Aut(\bar{\eta})= \{\pm I\} \subset \pi_1(\Ell_{\bar{k}}).\]

Let $\Tate/k((q))$ denote the base-change of the Tate curve to $k((q))$. There is a unique isomorphism 
\[k[[1/j]] \lisim k[[q]], \; j \longmapsto j(q):=j(\Tate/k((q)))=\tfrac{1}{q}+744+\cdots,\]
which is compatible with the adic topologies (cf. \cite[(8.11.3)]{katzArithmeticModuliElliptic1985}). Via the induced isomorphism
\[\spec(k((q))) \lisim \left(\A^{1}_{k}\right)^\wedge_\infty ,\] 
the Tate curve gives rise to a section of $j:\Ell_{k,\infty}^\wedge\to \left(\A^{1}_{k}\right)^\wedge_\infty$. By functoriality of $\pi_1$ and Remark \ref{rem:completions}, we obtain the following. 
\begin{prop}\label{prop:Tatecurvesplitting}
	The Tate curve induces a splitting of the exact sequence
	\[1 \to \{\pm I\}\to \pi_1(\Ell_{k,\infty})\to \pi_1(\A^{1}_{k,\infty})\to 1,\]
	such that \[\pi_1(\Ell_{k,\infty})\iso  \pi_1(\A^{1}_{k,\infty})\times \{\pm I\}.\]
\end{prop}
Consider the Henselian discrete valued field $F=k((q))$ equipped with the $q$-adic valuation, whose residue field is $k$. Let  $\overline{F}$ denote the algebraic closure and $F^{\text{nr}}\subset \overline{F}$ the maximal unramified extension. Associated with $F$, we have the short exact sequence 
\[1 \to \widehat{\Z}(1)=\Gal_{ F^{\text{nr}}} \to \Gal_{k((q))} \to \Gal_k \to 1,\]
which we denote $\pi_1(k((q))/k)$. Let $\pi= (q^{\frac{1}{n}})_{n \in \N} \in \overline{F}$ be a compatible system of roots of $q$, then
\[ F^{\text{nr}}\cap F(\pi)=F \quad \text{and} \quad \overline{F}=F^{\text{nr}}F(\pi).\]
Therefore $\pi$ determines a section
\[t_\pi:\Gal_{k}= \Gal(F^{\text{nr}}/F)=\Gal(F^{\text{nr}}F(\pi)/F(\pi)) \subset\Gal_{F}\] of the exact sequence $\pi_1(k((q))/k)$. Let $\mathscr{S}_{\pi_1(k((q))/k)}$ be the $\widehat{\Z}(1)$-conjugacy classes of sections of $\pi_1(k((q))/k)$. We identify \[ \mathscr{S}_{\pi_1(k((q))/k)}= \H^1(k,\widehat{\Z}(1))\]
as follows. Given $t\in \mathscr{S}_{\pi_1(k((q))/k)}$, we consider the difference cocycle 
\[\delta(t,t_\pi)=tt_\pi^{-1} \in \mathrm{Z}^1(\Gal_{k},\widehat{\Z}(1)).\]
Replacing $\pi$ with another compatible system, modifies $\delta(t,t_\pi)$ by a coboundary. Moreover the cohomology class satisfies
\[[\delta(t,t_\pi)]=t^*[q] \in \H^1(k,\widehat{\Z}(1)),\]
where $t^*[q]$ is the pullback of the class $[q] \in \H^1(k((q)),\widehat{\Z}(1))$ along $t$. We call 
\[ q(t):=t^*[q] \in  \H^1(k,\widehat{\Z}(1))\]
the \textbf{$q$-parameter of $t$}.

\vspace{\baselineskip}

The Tate curve defines a morphism $\Tate\colon \spec(k((q)))\to \Ell_{k,\infty}^\wedge$ with induced map on sections $\Tate\colon \mathscr{S}_{\pi_1(k((q))/k)}\to \mathscr{S}_{\pi_1(\Ell_{k,\infty}/k)}$. As an immediate consequence of Proposition \ref{prop:Tatecurvesplitting}, we have the following.

\begin{prop}\label{prop:classificationofcuspidals}
	The morphism  $\Tate\colon \spec(k((q)))\to \Ell_{k,\infty}^\wedge$ induces a bijection
\begin{align*}
\mathscr{S}_{\pi_1(k((q))/k)} \times \H^1(k,\{\pm 1\})&\lisim	\mathscr{S}_{\pi_1(\Ell_{k,\infty}/k)}\\ 
	 (t,\epsilon)&\longmapsto \Tate(t)\otimes \epsilon.
	\end{align*}
\end{prop}
Thus, for every $s \in  \mathscr{S}_{\pi_1(\Ell_{k,\infty}/k)}$ there is a unique section $t \in \mathscr{S}_{\pi_1(k((q))/k)} $ and a unique quadratic character $\epsilon_{s}$, such that
\[s=\Tate(t)\otimes \epsilon_{s}.\] 
We call $q(s):=q(t) \in \H^1(k,\widehat{\Z}(1))$ the \textbf{$q$-parameter} of $s$.

\begin{prop}\label{prop:ramificationindexforcuspidal}
Let $V/k$ be a geometrically connected curve with smooth completion $X$ and $y \in (X \setminus V)(k)$ a $k$-rational cusp. Let $f\colon V\to \Ell_{k}$ be a finite, surjective $k$-morphism and let $e=e_{y/\infty}$ be the ramification index with respect to $j\circ f$. Choose a uniformizer $u \in \widehat{\o}_{X,y}$. For any cuspidal section $s \in \mathscr{S}_{\pi_1(V_y/k)}$, we have
\[q(f(s))\in [r]+eH^1(k,\widehat{\Z}(1)),\]
where $[r]$ is the Kummer class associated to the residue
\[r= \restr{\frac{1/(j\circ f)}{u^e}}{y} \in k^\times.\]
\end{prop}
\bp 
Using the uniformizer $u$, we may identify 
\[V_y^\wedge =\spec(k((u)))\]
and write 
\[ 1/{j\circ f}= r \cdot u^e \cdot v\]
for some $v \in 1+u k[[u]]$. Since $k$ has characteristic zero, there exists $w\in 1+u k[[u]]$ such that $v=w^e$. Hence
\[
f^*([1/j])=[1/(j\circ f)]=\pr^*([r])+e[uw] \in \H^1(V_y^\wedge,\widehat{\Z}(1)),
\]
where $\pr\colon V_y^\wedge \to \spec(k)$ is the structure map.
Thus
\[q(f(s))=f(s)^*([1/j])=s^*(f^*([1/j]))=  [r]+e\cdot s^*([uw]) \in [r]+e\H^1(k,\widehat{\Z}(1)).\qedhere\]
\ep 

\subsection{Attaching Galois representations to sections}
Let $k$ be of characteristic zero. An elliptic curve $E_0/\Omega$ over an algebraically closed field extension $\Omega/k$ defines a geometric point $\bar{x} \in \Ell_{k}$.  The Tate module of the universal elliptic curve over $\Ell_k$ gives rise to a monodromy representation 
\[\rho\colon \pi_1(\Ell_k,\bar{x}) \longrightarrow \GL(TE_0)\iso \GL_2(\widehat{\Z}),
\]
where $TE_0= \prod_\l T_\l E_0$ is the Tate module of $E_0$. Since the Weil-pairing is compatible with base-change, it induces an isomorphism of $\pi_1(\Ell_k,\bar{x})$ modules 
\[\bigwedge^2TE_0 \lisim \widehat{\Z}(1).\] 
Hence the determinant is given by \[ \det( \rho) =\chi \circ  \pi_1(\pr),\] where $\pr \colon \Ell_{k} \to \spec(k)$ is the structure map and $\chi\colon \Gal_{k}\to \widehat{\Z}^\times$ is the cyclotomic character. For $N\geq 1$, we write 
\[ \overline{\rho}_N  \colon \pi_1(\Ell_k,\bar{x}) \longrightarrow \GL(E_0[N])\iso \GL_2(\Z/N\Z),
\]
for the induced representation on the $N$-torsion. For a section $s\in\mathscr S_{\pi_1(\Ell_k/k)}$, we define
\[
\rho_s:=\rho\circ s
\qquad\text{and}\qquad
\overline{\rho}_{s,N}:=\overline{\rho}_N\circ s.
\]
Then $\rho_s$ is a well-defined adelic Galois representation up to isomorphism, satisfying
\[\det(\rho_{s})=\chi\circ \pi_1( \pr)\circ s=\chi.\]
Writing $\overline{\chi}_N \colon \Gal_k\to (\Z/N\Z)^\times$ for the cyclotomic character mod $N$, we have 
\[\det(\overline{\rho}_{s,N} ) = \overline{\chi}_N.\]
If $s=s_E$ is the rational section associated to an elliptic curve $E/k$, then $\rho_s$ and $\overline{\rho}_{s,N}$ coincide with the natural Galois representations on $TE$ and $E[N]$, respectively. We denote these by $\rho_E$ and $\overline{\rho}_{E,N}$. 
\begin{lem}\label{lem:quadratictwistcompatibility}
	Let $s \in \mathscr{S}_{\pi_1(\Ell_k/k)}$ and $\epsilon\colon \Gal_k \to \{\pm 1\}$ be a quadratic character. Then 
	\[\rho_{s\otimes \epsilon}\iso \rho_{s}\otimes \epsilon.\]
\end{lem}
\bp 
Let $\bar{x} \in \Ell_{k}$ be a geometric point, which corresponds to an elliptic curve $E_0/\Omega$ for some algebraically closed extension $\Omega/k$. Consider the monodromy representation associated to the Tate module of the universal elliptic curve
\[\rho\colon \pi_1(\Ell_{k},\bar{x})\longrightarrow \GL(TE_0).\]
By Lemma \ref{lem:sl2}, the element 
\[-I\in  \pi_1(\Ell_{\bar{k}},\bar{x}) \iso \widehat{\SL_2(\Z)}\] is  the unique element of order 2. The image of the automorphism $-1 \in \Aut(E_0)$ under the homomorphism
\[ \Aut(E_0)=\Aut(\bar{x}) \xrightarrow{\omega_{\bar{x}}} \pi_1(\Ell_{\bar{k}},\bar{x}) \subset \pi_1(\Ell_{k},\bar{x})\]
is an element of order 2, hence coincides with $-I$. Since $-1 \in \Aut(E_0)$ acts on $TE_0$ via multiplication by $-1$, it follows that 
\[\rho(-I)=-\id_{TE_0} \in \GL(TE_0),\] which implies the assertion.
\ep 
\begin{prop}\label{prop:cuspidalrep}
	Let $s \in \mathscr{S}_{\pi_1(\Ell_{k,\infty}/k)}$ with $q$-parameter $q(s) \in \H^1(k,\widehat{\Z}(1))$ and quadratic character $\epsilon_s\colon \Gal_k \to \{\pm 1\}$. Then
	\[\rho_{s}\iso \rho_{q(s)}\otimes \epsilon_s,\]
	where $\rho_{q(s)}=\begin{pmatrix}
		\chi & \ast \\
		0 & 1
	\end{pmatrix}$
	is the extension corresponding to $q(s) \in \H^1(k,\widehat{\Z}(1))$.
\end{prop}
\bp 
By Proposition \ref{prop:classificationofcuspidals}, we have $s=\Tate(t) \otimes \epsilon_s$ with $t \in  \mathscr{S}_{\pi_1(k((q))/k)}$. By Lemma \ref{lem:quadratictwistcompatibility}, it suffices to treat the case of $\epsilon_s=1$. For $N\geq 1$ let $\Tate[N]/k((q))$ denote the $N$-torsion group scheme of the Tate-curve over $k((q))$. By \cite[Section VII]{deligneSchemasModulesCourbes1973}, we have a canonical extension

\[ 0\to \mu_N \to  \Tate[N]  \to \Z/N \Z \to 0 \]
corresponding to the Kummer class $[q] \in \H^1(k((q)),\mu_N)$. Thus, considered as extensions of $\Z/N\Z$ by $\mu_N$, we have
\[t^*\Tate[N]=t^*[q]=q(t)=q(s) \in \H^1(k,\mu_N). \qedhere\]
\ep 
\section{Selmer sections of the moduli stack of elliptic curves over \texorpdfstring{$\Q$}{Q}}\label{sec:selmersec}
\subsection{Finite descent obstruction for modular curves} For $N\geq 1$, let $\Y_1(N),\Y(\mu_N)$ and $\Y_0(N)$ denote the finite \'etale covers of $\Ell_{\Q}$ parametrizing elliptic curves equipped, respectively, with a $\Z/N\Z$-subgroup, a $\mu_N$-subgroup and a cyclic subgroup of order $N$ in the $N$-torsion. The respective coarse moduli spaces $Y_1(N),Y(\mu_N)$ and $Y_0(N)$ are geometrically connected curves over $\Q$ with smooth completions $X_1(N),X(\mu_N)$ and $X_0(N)$. By \cite[Chapter V, Proposition 1.4]{deligneSchemasModulesCourbes1973}, there is an isomorphism $\mathcal{Y}_1(N)\iso \mathcal{Y}(\mu_N)$, which induces an isomorphism \[X_1(N)\iso X(\mu_N).\]

The next proposition is due to Stoll, see \cite[Corollary 8.8.]{stollFiniteDescentObstructions2007}. 
\begin{prop}\label{prop:stollmodular}
	 If $X \in \{X_1(N),X(\mu_N),X_0(N)\}$ has positive genus, then 
	\[X(\A_\Q)^{\mathrm{f\text{-}cov}}=X(\Q).\]
\end{prop}
Stoll's statement does not include $X(\mu_N)$, but this follows immediately from the isomorphism $X(\mu_N)\iso X_1(N)$.

Recall that the modular curves $X_0(N)$ (resp.\ $X_1(N)\iso X(\mu_N))$ have genus zero if and only if $N\in S_0$ (resp.\ $N\in S_1$), where 
\[ S_0:= \{1,2,\dots,10\} \cup \{12,13,16,18,25\} \quad \text{and} \quad S_1:= \{1,2,\dots,10\} \cup \{12\}.\]
Let $G$ be a group and $M$ a $(\mathbb{Z}/N\mathbb{Z})[G]$-module, which is finite and free as a $\mathbb{Z}/N\mathbb{Z}$-module. We call $M$  \textbf{reducible}, if there exists a proper $(\mathbb{Z}/N\mathbb{Z})[G]$-submodule of $M$, which is itself free over $\mathbb{Z}/N\mathbb{Z}$.
\begin{thm}\label{thm:finitedescent}
	Let $s\in \mathscr{S}_{\pi_1(\Ell_\Q/\Q)}$ be a Selmer section. Assume that one of the following conditions is satisfied.
	\be[label=(\alph*)]
	\item The representation $\overline{\rho}_{s,N}$ is reducible and $N\notin S_0$. 
	\item The representation $\overline{\rho}_{s,N}$ has $\Z/N\Z$ or $\mu_N$ as a sub-representation and $N \notin S_1$.  	
	\ee
	Then $s$ is rational or cuspidal. 
\end{thm}
\bp 
Assume (a). After choosing base-points, consider the monodromy representation on the $N$-torsion of the universal elliptic curve\[
\overline{\rho}_N \colon \pi_1(\Ell_\Q)\longrightarrow \GL_2(\Z/N\Z).
\]
By definition, the subgroup $\pi_1(\mathcal Y_0(N))\subset \pi_1(\Ell_{\Q})$ consists of elements $g \in \pi_1(\Ell_{\Q})$ preserving a fixed cyclic subgroup $ C\subset (\Z/N\Z)^2$ of order $N$. Since $\overline{\rho}_{s,N}$ is reducible, its image preserves a cyclic subgroup of order $N$. After conjugating, we conclude that
\[
s(\Gal_{\Q})\subseteq \pi_1(\mathcal Y_0(N)).
\]
Thus $s$ lifts to a section $\tilde{s} \in \mathscr{S}_{\pi_1(\mathcal{Y}_0(N)/\Q)}$. By Lemma \ref{lem:goingupanddown}, the section $\tilde{s}$ is Selmer again. Since $X_0(N)$ has positive genus, Proposition \ref{prop:stollmodular} implies
\[X_0(N)(\A_\Q)^{\mathrm{f\text{-}cov}}=X_0(N)(\Q).\]
Applying Proposition \ref{prop:finitedescentmoduli} to the non-constant $\Q$-morphism
\[\mathcal{Y}_0(N) \xrightarrow{\text{coarse moduli}} Y_0(N) \hookrightarrow X_0(N),\]
we conclude that $\tilde{s}$ is cuspidal or rational. Hence $s$ is cuspidal or rational. The proof of (b) is analogous.
\ep

\subsection{$\l$-cuspidal sections}
Let $\l$ be a prime. Given a finite-dimensional representation
\[
\overline{\rho}\colon G\longrightarrow \GL_n(\F_\ell)
\]
of a profinite group $G$, let $\overline{\rho}^{\ss}$ denote its semisimplification. For every closed subgroup $H\subset G$, we write $\left.\overline{\rho}\right|_{H}^{\ss}$ for the semisimplification of the restriction $\restr{\overline{\rho}}{H}$. Furthermore, let
\[
\overline{\chi}_\ell\colon \Gal_{\Q}\longrightarrow \F_\ell^\times
\]
denote the reduction of the $\ell$-adic cyclotomic character
\[
\chi_\ell\colon \Gal_\Q \longrightarrow \Z_\ell^\times.
\]
By abuse of notation, we use the symbol $\overline{\chi}_\ell$ also for its restriction to any decomposition or inertia subgroup. Finally, for a prime $p$, let $I_p\subset\Gal_{\Q_p}$ denote the inertia subgroup.

We now introduce the notion of an \textit{$\l$-cuspidal} section in $\mathscr{S}_{\pi_1(\Ell_\Q/\Q)}$. This will allow us to determine whether $s \in \mathscr{S}_{\pi_1(\Ell_\Q/\Q)}$ is cuspidal in terms of the ramification of $\overline{\rho}_{s,\l}$. 
\begin{defi}
	Let $\l$ be a prime.
	\be 
	\item  We say that $s \in \mathscr{S}_{\pi_1(\Ell_{\Q_p}/\Q_p)}$ is \textbf{$\l$-cuspidal}, if 
	there exists a quadratic character $\epsilon\colon  \Gal_{\Q_p}\to \{\pm 1\}$, such that 
	\[\left.\overline{\rho}_{s,\l}\right|_{I_p}^{\ss} = \epsilon \oplus \epsilon\overline{\chi}_\l.\]
	\item  A section $s\in \mathscr{S}_{\pi_1(\Ell_\Q/\Q)}$ is called \textbf{$\l$-cuspidal}, if  $s_{\Q_p}$ is $\l$-cuspidal for all primes $p$ and $\overline{\rho}_{s,\l}$ is reducible.
	\ee 
\end{defi}
By Proposition \ref{prop:cuspidalrep}, a cuspidal section is $\l$-cuspidal.
\begin{lem}\label{lem:quadraticcharacters}
	Let $\Sigma$ be a finite set of rational primes and $\epsilon_p\colon \Gal_{\Q_p} \to \{\pm 1\}$ a quadratic character for each $p \in \Sigma$. Then there is a quadratic character $\epsilon\colon \Gal_\Q \to \{\pm 1\}$ unramified outside $\Sigma$, such that 
	\[ \restr{\epsilon}{I_p}= \restr{\epsilon_p}{I_p}\quad \text{	for all $p \in \Sigma$.}\]
\end{lem}
\bp 
By class field theory we have a commutative diagram 
\[\begin{tikzcd}
	{\prod_pI_{\Q_p^\ab/\Q_p}} && {\Gal_\Q^\ab} \\
	{\prod_p\Z_p^\times} && {\widehat{\Z}^\times,}
	\arrow["\sim", from=1-1, to=1-3]
	\arrow["{\chi^{-1},\sim}", from=1-3, to=2-3]
	\arrow["\sim"', from=2-1, to=1-1]
	\arrow["\sim", from=2-1, to=2-3]
\end{tikzcd}\]
where the left, vertical isomorphisms are induced by the local reciprocity homomorphism $\Q_p^\times \to \Gal_{\Q_p}^\ab$ and $\chi^{-1}$ is the reciprocal of the cyclotomic character. The characters $\epsilon_p$ restrict to characters $\tilde{\epsilon}_p \colon I_{\Q_p^\ab/\Q_p} \to \{\pm 1\}$. Because $\Sigma$ is finite, the product 
\[ \epsilon \colon \Gal_\Q^\ab \iso  \prod_pI_{\Q_p^\ab/\Q_p} \xrightarrow{\prod_{p \in \Sigma}\tilde{\epsilon}_p} \{\pm 1\}\]
is well defined and satisfies the desired properties.
\ep 
\begin{prop}\label{prop:lcuspidalimplieslift}
	A section $s\in \mathscr{S}_{\pi_1(\Ell_\Q/\Q)}$ is $\l$-cuspidal if and only if there is a  quadratic character $\epsilon\colon \Gal_\Q \to \{\pm 1\}$, such that $s\otimes \epsilon$ lifts to a section of $\mathcal{Y}_1(\l)$ or $\mathcal{Y}(\mu_\l)$. 
\end{prop}
\bp 
Let $s \in \mathscr{S}_{\pi_1(\Ell_\Q/\Q)}$ be an $\l$-cuspidal section. For each prime $p$, let $\epsilon_p\colon \Gal_{\Q_p} \to \{\pm 1\}$ be a quadratic character, such that
\[ \left.\overline{\rho}_{s,\l}\right|_{I_p}^{\ss}=\left.\overline{\rho}_{s_{\Q_p},\l}\right|_{I_p}^{\ss}= \epsilon_p \oplus \epsilon_p\overline{\chi}_\l.\]
Let $\Sigma$ be the set where $\overline{\rho}:=\overline{\rho}_{s,\l}$ is ramified. In particular, the characters $\epsilon_p$ are unramified for $p\notin \Sigma$. By Lemma \ref{lem:quadraticcharacters}, there is a quadratic character $\epsilon\colon \Gal_\Q\to \{\pm 1\}$ unramified outside $\Sigma$, such that 
\[\restr{\epsilon}{I_p}=\restr{\epsilon_p}{I_p}\]
for all $p$.
Replacing $s$ with $s\otimes \epsilon$ we may assume that $\epsilon_p$ is unramified for all $p$. Since $\overline{\rho}$ is reducible, there is a one dimensional sub-representation $\omega$. Then 
\[\overline{\rho}^{\ss}=\omega^{-1} \overline{\chi}_\l \oplus \omega.\]
By assumption we have
\[\restr{\overline{\rho}^{\ss}}{I_p}= 1 \oplus \overline{\chi}_\l \quad \text{for all $p$.}   \]
We conclude that $\omega$ is unramified outside $\l$. At $p=\l$ either $\omega$ or $\omega^{-1} \overline{\chi}_\l$ is unramified. By Minkowski's Theorem, $\Q$ has no non-trivial unramified extensions, which implies that either $\omega=1$ or $\omega= \overline{\chi}_\l$. Therefore $s$ lifts to $\mathcal{Y}_1(\l)$ or to $\mathcal{Y}(\mu_\l)$. The converse is clear.
\ep 

\begin{cor}
	Let  $s\in \mathscr{S}_{\pi_1(\Ell_\Q/\Q)}$ be a Selmer section and $\l\geq 11$ a prime.  Then $s$ is cuspidal if and only if $s$ is $\l$-cuspidal.
\end{cor}
\bp 
A cuspidal section is  $\l$-cuspidal by Proposition \ref{prop:cuspidalrep}. Let $s$ be a $\l$-cuspidal Selmer section. By Proposition \ref{prop:lcuspidalimplieslift}, we may replace $s$ with a quadratic twist and assume that $s$ lifts to $\mathcal{Y}_1(\l)$ or $\mathcal{Y}(\mu_\l)$. By Theorem \ref{thm:finitedescent}, we conclude that $s$ is cuspidal or rational. By Mazur's Torsion Theorem \cite[Theorem (8)]{mazurModularCurvesEisenstein1977}  the stacks $\mathcal{Y}_1(\l)$ and  $\mathcal{Y}(\mu_\l)$ have no $\Q$-points, implying that $s$ is cuspidal.
\ep

\begin{lem}\label{l-cuspidalellipticcurves}
	Let $s\in \mathscr{S}_{\pi_1(\Ell_{\Q_p}/\Q_p)}$ be a rational section associated to an elliptic curve $E/\Q_p$ and $\l$ a prime.  
	\be 
	\item If $E/\Q_p$ has potential multiplicative reduction, then $s$ is $\l$-cuspidal.
	\item Assume that $E/\Q_p$ has good reduction. 
	\begin{itemize}
		\item  If $\l\neq p$, then $s$ is $\l$-cuspidal. 
		\item If $\l=p$, then $s$ is $\l$-cuspidal if and only if $E/\Q_p$ has good ordinary reduction.
	\end{itemize}
	\ee 
\end{lem}
\bp 
(1) If $E/\Q_p$ has potential multiplicative reduction, then $E$ is a quadratic twist of a Tate curve and we have 
\[\overline{\rho}_{s,\l}\iso    \begin{pmatrix} \overline{\chi}_\l & \ast \\
	0 & 1 
\end{pmatrix}\otimes \epsilon\]
for some quadratic character $\epsilon \colon \Gal_{\Q_p} \to \{\pm 1\}$. Therefore $s$ is $\l$-cuspidal.

(2) Suppose that $E/\Q_p$ has good reduction. If $\l\neq p$, then $\overline{\rho}_{s,\l}$ is unramified and $s$ is clearly $\l$-cuspidal. Assume that $\l=p$.  If $E$ has ordinary reduction, then $ \left.\overline{\rho}_{s,\l}\right|_{I_\l}$ has a trivial one-dimensional quotient by \cite[Proposition 11]{serreProprietesGaloisiennesPoints1971}. Hence $s$ is $\l$-cuspidal.  If $E$ has supersingular reduction, the image of inertia $\overline{\rho}_{s,\l}(I_\l) \subset \GL_2(\F_\l)$ is a non-split Cartan subgroup by \cite[Proposition 12]{serreProprietesGaloisiennesPoints1971}, which is irreducible. Therefore $s$ cannot be $\l$-cuspidal.
\ep

\subsection{Selmer sections of the Fermat curve}
Let $\l$ be a prime. The \textbf{Fermat curve} $V_\l/\Q$ is the affine hyperbolic curve defined by the homogeneous equations
\[X^\l + Y^\l+ Z^\l = 0\quad  \text{and} \quad XYZ \neq 0
\]
in $\P^2_\Q$. To ease notation, we will write $\mathscr{S}_{\pi_1(V_\l/\Q_p)}:=\mathscr{S}_{\pi_1(V_{\l}\times_\Q \Q_p/\Q_p)}$ for primes $p$.
 We view $V_\l$ as a finite étale cover of $\mathbb{P}^1_\mathbb{Q} \setminus \{0,1,\infty\}$ via the map
\[
	V_\l \longrightarrow \mathbb{P}^1_\mathbb{Q} \setminus \{0,1,\infty\}, \;
	[X:Y:Z] \longmapsto -X^\l/Z^\l.
\]
The \textbf{Legendre family} of elliptic curves  $\mathcal{E}_{\mathrm{Leg}}$ over $\P^1_\Q\setminus \{0,1,\infty\}$ is the family of elliptic curves with affine equation 
\[y^2=x(x-1)(x-\lambda) \quad \lambda \in  \P^1_\Q\setminus \{0,1,\infty\}.\]
The classifying morphism 
\[\P^1_\Q\setminus \{0,1,\infty\}  \xrightarrow{\mathcal{E}_{\mathrm{Leg}}}  \Ell_{\Q} \]
is finite \'etale and realizes $V_\l$ as a finite \'etale cover of $\Ell_{\Q}$ via the composition 
\[V_\l \longrightarrow\P^1_\Q\setminus \{0,1,\infty\} \xrightarrow{\mathcal{E}_{\mathrm{Leg}}} \Ell_{\Q}. \]
Using this finite \'etale cover, we attach a Galois representation $\rho_{s,\l}$ to a section $s \in \mathscr{S}_{\pi_1(V_\l/\Q)}$ (resp.\ $\mathscr{S}_{\pi_1(V_\l/\Q_p)})$, by considering its image in $\mathscr{S}_{\pi_1(\Ell_{\Q}/\Q)}$ (resp.\ $\mathscr{S}_{\pi_1(\Ell_{\Q_p}/\Q_p)})$. In particular, the property of $s$ being cuspidal (resp.\ rational) is preserved by (quadratic) twisting by Lemma~\ref{lem:beingcuspidaldescends} and going up or down along finite \'etale covers by Lemma \ref{lem:goingupanddown}.

\vspace{\baselineskip}

Let $\overline{\rho} \colon \Gal_{\mathbb Q_p} \to \GL_2(\overline{\mathbb F}_\ell)$ be a two-dimensional residual Galois representation. Recall that if $\ell \neq p$, one associates to $\overline{\rho}$ a \textbf{local conductor} $n_p(\overline{\rho})$, whereas if $\ell = p$, one associates a \textbf{weight} $k(\overline{\rho})$ (see \cite[§1-2]{serreRepresentationsModulairesDegre1987}).
 In both cases, these invariants depend only on the restriction to the inertia group $\restr{\overline{\rho}}{I_p}$.

\begin{defi}
	Let $\l>2$ be a prime.
	\be 
	\item We say that a section $s \in \mathscr{S}_{\pi_1(\Ell_{\Q_p}/\Q_p)}$ is \textbf{$\l$-Fermat}, if there exists a quadratic character $\epsilon: \Gal_{\Q_p}\to \{\pm 1\}$ such that the following holds.
	\begin{itemize}
		\item  If $p=2$, the local conductor satisfies 
		\[n_2(\restr{\overline{\rho}_{s\otimes \epsilon,\l}}{I_2})\leq 4.\] 
		\item If $p\notin  \{2,\l\}$, the representation $\overline{\rho}_{s\otimes \epsilon,\l}$ is unramified. 
		\item If $p= \l$, the representation $\left.\overline{\rho}_{s\otimes \epsilon,\l}\right|_{I_\l} $ is of weight 2.
	\end{itemize}  
	\item A section $s  \in \mathscr{S}_{\pi_1(\Ell_{\Q}/\Q)}$ is called \textbf{$\l$-Fermat}, if $s_{\Q_p}$ is $\l$-Fermat for all primes $p$ and $\overline{\rho}_{s,\l}$ is irreducible.
	\ee 
\end{defi}
\begin{prop}\label{prop:nolfermatsections}
	There are no $\l$-Fermat sections in $\mathscr{S}_{\pi_1(\Ell_{\Q}/\Q)}$ for $\l>2$.
\end{prop}
\bp 
Let $s \in \mathscr{S}_{\pi_1(\Ell_{\Q}/\Q)}$ be an $\l$-Fermat section and let $\overline{\rho}=\overline{\rho}_{s,\l}$. Since $\overline{\rho}$ is a 2-dimensional, odd and irreducible, $\overline{\rho}$ is absolutely irreducible. Indeed, if $c \in \Gal_\Q$ is a complex conjugation, then $\overline{\rho}(c)$ has order 2 with eigenvalues $-1$ and $1$. Hence any one-dimensional invariant subspace must be an eigenspace of $\overline{\rho}(c)$, thus already defined over $\F_\l$.

Let $\Sigma_0$ be the finite set of primes $p$, such that $p \notin \{2,\l\}$ and $\overline{\rho}$ is ramified at $p$. For each $p \in \Sigma :=\Sigma_0 \cup \{2,\l\}$, we find quadratic characters $\epsilon_p\colon \Gal_{\Q_p} \to \{\pm 1\}$, such that  
\begin{itemize}
	\item  at $p=2$, the local conductor satisfies
	\[n_2(\restr{\overline{\rho}_{s_{\Q_2}\otimes \epsilon_2,\l}}{I_2})\leq 4,\] 
	\item for $p\in \Sigma_0$, the representation $\overline{\rho}_{s_{\Q_p}\otimes \epsilon_p,\l}$ is unramified and
	\item for $p= \l$, the representation $\left.\overline{\rho}_{s_{\Q_\l}\otimes \epsilon_\l,\l}\right|_{I_\l} $ is of weight 2.
\end{itemize}
 By Lemma \ref{lem:quadraticcharacters}, there is a quadratic character $\epsilon \colon \Gal_{\Q} \to \{\pm 1\}$ unramified outside $\Sigma$, such that for all $p \in \Sigma$, we have
\[ \restr{\epsilon}{I_p}=\restr{\epsilon_p}{I_p}.\]
Replacing $s$ with $s\otimes \epsilon$, we find that $\overline{\rho}$ has weight 2 at $\l$, is unramified outside $\{2,\l\}$ and the level $N(\overline{\rho})$ divides $2^4=16$. Applying Serre's Conjecture, proven by Khare and Wintenberger \cite{khareSerresModularityConjecture2009}, we conclude that $\overline{\rho}$ arises from a normalized eigenform $f$ of weight $2$ and level $N(\overline{\rho})$. The modular curve $X_0(N)$ has genus 0 for all $N \mid 16$, contradicting the existence of $f$.
\ep 
We say that a section $s \in \mathscr{S}_{\pi_1(V_\l/\Q_p)}$ is \textbf{$\l$-cuspidal} (resp.\ \textbf{$\l$-Fermat}), if the image of $s$ in $ \mathscr{S}_{\pi_1(\Ell_{\Q_p}/\Q_p)}$ is  $\l$-cuspidal (resp.\ $\l$-Fermat).
\begin{prop}\label{prop:cuspidalsectionfermat}
	Every cuspidal section $s \in \mathscr{S}_{\pi_1(V_\l/\Q_p)}$ is $\l$-Fermat for $\l>2$.
\end{prop}
\bp 
Let $s \in \mathscr{S}_{\pi_1(V_\l/\Q_p)}$ be a cuspidal section associated to a $\Q_p$-rational cusp $y$ of $V_\l$. Recall that by definition $s$ is $\l$-Fermat if $s_0:=f(s) \in \mathscr{S}_{\pi_1(\Ell_{\Q_p,\infty}/\Q_p)}$ is $\l$-Fermat, where
\[f\colon V_\l \longrightarrow \P^1_{\Q}\setminus \{0,1,\infty\}\xrightarrow{\mathcal{E}_{\mathrm{Leg}}} \Ell_{\Q}\] 
is the composition of the previously fixed finite \'etale cover $V_\l \to \P_\Q^1\setminus \{0,1,\infty\}$ and the classifying map of the Legendre family.

We begin by analyzing the $q$-parameter of $s_0$. Let $j\colon \Ell_{\Q} \to \A^1_\Q$ be the $j$-invariant. The map $j\circ f$ is given by the formula 
\[ j(\lambda)=2^8\frac{(\lambda^2-\lambda+1)^3}{\lambda^2(\lambda-1)^2}, \quad \lambda = -X^\l/Z^\l. \]
The finite \'etale cover
\[
V_\l\longrightarrow \P^1_{\Q}\setminus\{0,1,\infty\},
\qquad
[X:Y:Z]\longmapsto -X^\l/Z^\l,
\]
extends to a map of smooth compactifications with ramification index $\l$ at every cusp of $V_\l$. The composition
\[
\P^1_\Q\setminus \{0,1,\infty \}\xrightarrow{\mathcal{E}_{\mathrm{Leg}}} \Ell_{\Q}\xrightarrow{j} \A^1_{\Q}
\]
has ramification index $2$ at each of the cusps $0,1,\infty$. Thus the ramification index of $j\circ f$ at $y$ is
\[
e_{y/\infty}=2\l.
\]
Let
\[
X_\l\colon X^\l+Y^\l+Z^\l=0
\]
be the smooth completion of $V_\l$. Depending on whether $X,Y$ or $Z$ vanishes on $y$, we choose the uniformizer $u \in \widehat{\o}_{X_\l,y}$ to be 
\[X/Z, \quad Y/Z, \quad Z/X. \]
In any case, one of
\[\lambda=-X^\l/Z^\l, \quad 1-\lambda=-Y^\l/Z^\l , \quad 1/\lambda=-Z^\l/X^\l\]
is equal to $-u^\l$. Since 
\[j(\lambda)=j(1-\lambda)=j(1/\lambda),\] we have 
\[\frac{1}{j \circ f} = \frac{1}{j(-u^\l)} =2^{-8}u^{2\l}\frac{(1+u^\l)^2}{(1+u^\l+u^{2\l})^3}.\]
 By Proposition \ref{prop:ramificationindexforcuspidal}, the $q$-parameter satisfies
\[q(s_0) \in [2^{-8}] +2\l \H^1(\Q_p,\widehat{\Z}(1)).\]
In particular, its image in $\H^1(\Q_p,\F_\l(1))$ is the Kummer class $[2^{-8}]$.
Applying Proposition \ref{prop:cuspidalrep} and reducing modulo $\l$, we obtain
\[\overline{\rho}_{s_0,\l}\iso \overline{\rho}\otimes \epsilon.\] Here
\[\overline{\rho}\iso  \begin{pmatrix}
	\overline{\chi}_\l & \ast \\
	0 & 1
\end{pmatrix}\] is the representation associated to the Kummer class $[2^{-8}] \in \H^1(\Q_p,\F_\l(1))$. The character $\epsilon=\epsilon_{s_0}$ is the quadratic character associated to $s_0$ as in Proposition \ref{prop:classificationofcuspidals}.

The property of being $\l$-Fermat is invariant under quadratic twists. Therefore it suffices to prove that the section 
\[t:=s_0\otimes \epsilon \in \mathscr{S}_{\pi_1(\Ell_{\Q_p}/\Q_p)}\]
is $\l$-Fermat. We have
\[ \overline{\rho}_{t,\l} \iso \overline{\rho}_{s_0,\l}\otimes \epsilon \iso \overline{\rho}.\]

Case 1: $p=2$. Since $\l>2$, the representation $\restr{\overline{\rho}}{I_2}$ is unipotent and hence tamely ramified. Thus the local conductor satisfies
\[n_2(\overline{\rho})\leq 2 \]
and hence $t$ is $\l$-Fermat.

Case 2: $p \notin \{2,\l\}$. In this case $[2^{-8}] \in \H^1(\Q_p,\F_\l(1))$ is an unramified class. Hence
\[ \restr{\overline{\rho}}{I_p}\iso 1\oplus \overline{\chi}_\l.\] 
Since $p\neq \l$, the cyclotomic character $\overline{\chi}_\l$ is unramified. Thus $\overline{\rho}$ is unramified and therefore $t$ is $\l$-Fermat.

Case 3: $p=\l$. Since $\l>2$, we have
 \[ v_\l(2^{-8})=0.\]
Therefore $\overline{\rho}$ is peu ramifi\'ee and hence has weight 2 by \cite[Proposition 3]{serreRepresentationsModulairesDegre1987}. Hence $t$ is $\l$-Fermat.
\ep 
\begin{lem}\label{lem:valuationofj-inv}
	Let $\lambda \in \Q_p \setminus \{0,1\}$ and $E/\Q_p$ be the elliptic curve given by the equation  
	\[E\colon y^2=x(x-1)(x-\lambda).\]
If $v_p(\lambda)\neq 0$, then 
	\[v_p(j(E))=8 v_p(2)-2|v_p(\lambda)|.
	\]
\end{lem}
\bp 
We have 
\[j(E)=2^8\frac{(\lambda^2-\lambda+1)^3}{\lambda^2(\lambda-1)^2}.\]
We may assume that $v_p(\lambda)>0$, as replacing $\lambda$ by $1/\lambda$ changes $E$ by a quadratic twist and hence preserves the $j$-invariant. 
Since $v_p(\lambda)>0$, we have 
\[v_p(\lambda -1)=v_p(\lambda^2-\lambda+1)=0.\]
We conclude that
\[v_p(j(E))=8v_p(2)+3v_p(\lambda^2-\lambda+1)-2v_p(\lambda)-2v_p(\lambda-1)=8v_p(2)-2v_p(\lambda).\qedhere\]
\ep 
\begin{lem}\label{lem:reductionfermat}
	Let $E/\Q_p$ be an elliptic curve in the image of  $V_\l(\Q_p)\to \Ell(\Q_p)$. 
	\be 
		\item If $p=2$ and $\l\geq 5$, then after a quadratic twist $E$ is a Tate curve.
	\item If $p>2$, then after a quadratic twist either $E$ is a Tate curve with $\l \mid v_p(j(E))$ or $E$ has good reduction.
	\ee 
\end{lem}
\bp 
Let $z \in V_\l(\Q_p)$. In affine coordinates $z$ corresponds to tuple a $(a,b) \in \Q_p^2$ with $a^\l+b^\l+ 1=0$ and $ab\neq 0$. Under the finite \'etale cover
\[V_\l \longrightarrow \P_{\Q_p}^1\setminus \{0,1,\infty\}\to \Ell_{\Q_p},\]
the point $z$ is mapped to the elliptic curve $E/\Q_p$ with equation of the form
\[E \colon y^2=x(x-1)(x-\lambda),\]
where $\lambda =-a^\l$.

Case 1: $v_p(a)\neq 0$.  By Lemma \ref{lem:valuationofj-inv}, we have 
	\[v_p(j(E))=8 v_p(2)-2|v_p(\lambda)|=8v_p(2)-2\l|v_p(a)|<0\]
	for $p>2$ or $p=2$ and $\l\geq 5$. Hence $E$ is a Tate curve up to a quadratic twist. If $p>2$, then in particular $\l $ divides  $v_p(j(E))$.
	
Case 2: $v_p(a)=0$.  If $\lambda  \not\equiv 0,1 \bmod p$, then necessarily $p>2$ and the equation
\[y^2=x(x-1)(x-\lambda)\]
has good reduction modulo $p$. Thus, $E$ has good reduction and we are done.  We now assume that
\[ -a^\l \equiv  \lambda \equiv 0,1 \bmod  p.\]
Since $v_p(a)=0$, we conclude $\lambda \equiv 1 \bmod  p$. The equation 
\[a^\l +b^\l +1=0\]
implies that  $b \in \Z_p$ with $v_p(b)>0$. Replacing $\lambda$ with $\lambda':=1-\lambda$ produces a quadratic twist $E'$ of $E$. Since 
\[v_p(\lambda')=v_p(1-\lambda)=v_p(1+a^\l)=v_p(b^\l)>0,\]
we apply Lemma \ref{lem:valuationofj-inv} to $E'$ and compute
\[v_p(j(E))=v_p(j(E'))=8v_p(2)-2|v_p(\lambda')|=8v_p(2)-2\l|v_p(b)|<0.\]
Therefore $E$ is a quadratic twist of a Tate curve. If $p>2$, then $\l \mid v_p(j(E))$. 
\ep 
\begin{prop}\label{prop:rationalsectionsoffermat}
	Let $s \in\mathscr{S}_{\pi_1(V_\l/\Q_p)}$ be a rational section and $\l \geq 5$.
	\be 
	\item The section $s$ is $\l$-Fermat.
	\item If $\overline{\rho}_{s,\l}$ is reducible, then $s$ is $\l$-cuspidal.
	\ee 
\end{prop}
\bp 
Let $E/\Q_p$ be the elliptic curve associated to $s$. By Lemma \ref{lem:reductionfermat}, we may replace $E$ with a quadratic twist and assume that 
\begin{itemize}
	\item $p=2$ and $E$ is a Tate curve or
	\item $p>2$ and either  $E$ is a Tate curve with $\l \mid v_p(j(E))$ or $E$ has good reduction.
\end{itemize}
(1) We split up the proof in three cases.

Case 1: Assume that $p=2$ and $E$ is a Tate curve. Since $\l \neq 2$, the representation $\restr{\overline{\rho}_{E,\l}}{I_2}$ is unipotent. This means that $\overline{\rho}_{E,\l}(I_2)$ is an $\l$-group, so $\overline{\rho}_{E,\l}$ is tamely ramified. Hence the local conductor simplifies to\[n_2(\overline{\rho}_{E,\l} )= \dim_{\F_\l} E[\l]/E[\l]^{I_2}\leq 2.\]
It follows that $s$ is $\l$-Fermat.

Case 2: Assume that $p>2$ and $p\neq \l$. If $E$ has good reduction, then $\overline{\rho}_{E,\l}$ is unramified and in particular $s$ is $\l$-Fermat. Otherwise $E\iso E_q$ is a Tate curve with $q \in \Q_p^\times$ with $|q|<1$ and
\[\l \mid v_p(j(E))=v_p(j(E_q))=-v_p(q).\]
The class $[q] \in \im(\Q_p^\times\to \H^1(\Q_p,\mu_\l))$ describes the canonical extension 
\[0 \to \mu_\l  \to E_q[\l]\to \F_\l \to 0.\]
Hence $\l \mid v_p(q)$ implies that 
\[\restr{\overline{\rho}_{E_q,\l}}{I_p} \iso  \overline{\chi}_\l \oplus 1. \]
Therefore $\overline{\rho}_{s,\l}\iso \overline{\rho}_{E_q,\l}$ is unramified, which implies that $s$ is $\l$-Fermat. 

Case 3: Assume that $p>2$ and $p=\l$. Then $\overline{\rho}_{E,\l}$ has weight 2 by \cite[Proposition 5]{serreRepresentationsModulairesDegre1987}, therefore $s$ is $\l$-Fermat. 

(2) If $E$ is a Tate curve, then $s$ is $\l$-cuspidal, by Lemma \ref{l-cuspidalellipticcurves} (1). Assume that $E$ has good reduction. If $\l \neq p$, then $s$ is $\l$-cuspidal by Lemma \ref{l-cuspidalellipticcurves} (2). If $\l=p$, the reducibility of $\overline{\rho}_{s,\l}\iso \overline{\rho}_{E,\l}$ implies that $E$ has ordinary reduction. Hence, $s$ is $\l$-cuspidal by Lemma \ref{l-cuspidalellipticcurves} (2).
\ep 
\begin{thm}\label{thm:mainresult}
Let $\l\geq 7$ be a prime. Every Selmer section $s \in \mathscr{S}_{\pi_1(V_\l/\Q)}$ is cuspidal.
\end{thm}
\bp 
Let $s \in \mathscr{S}_{\pi_1(V_\l/\Q)}$ be a Selmer section. 

Case 1: Assume that $\overline{\rho}_{s,\l}$ is irreducible. For every prime $p$, the base-changed section $s_{\Q_p}$ is either cuspidal or rational. By Proposition \ref{prop:cuspidalsectionfermat} and Proposition \ref{prop:rationalsectionsoffermat} (1) the section $s_{\Q_p}$ is $\l$-Fermat in both cases. Therefore $s$ is $\l$-Fermat, which is a contradiction by Proposition \ref{prop:nolfermatsections}.

Case 2: Assume that $\overline{\rho}_{s,\l}$ is reducible. The monodromy representation on the 2-torsion of the Legendre family is trivial. Pulling back along $s$, we conclude that $\overline{\rho}_{s,2}$ is trivial. This implies that the representation $\overline{\rho}_{s,2\l}$ is reducible. Since $2\l >12$, the curve $X_0(2\l)$ has positive genus. Hence the section $s$ is cuspidal or rational by Theorem \ref{thm:finitedescent}.

Again, for every prime $p$, the section $s_{\Q_p}$ is cuspidal or rational. If $s_{\Q_p}$ is cuspidal, then $s_{\Q_p}$ is $\l$-cuspidal by Proposition \ref{prop:cuspidalrep}. Since $\overline{\rho}_{s,\l}$ is reducible, so is 
\[ \restr{\overline{\rho}_{s,\l}}{\Gal_{\Q_p}}\iso \overline{\rho}_{s_{\Q_p},\l}.\] 
Thus, if $s_{\Q_p}$ is rational, then $s_{\Q_p}$ is $\l$-cuspidal by Proposition \ref{prop:rationalsectionsoffermat} (2). We conclude that $s$ is $\l$-cuspidal.  By Proposition \ref{prop:lcuspidalimplieslift}, there is a quadratic character $\epsilon\colon \Gal_\Q \to \{\pm1\}$ such that $s \otimes \epsilon$ lifts to $\mathcal{Y}_1(\l)$ or $\mathcal{Y}(\mu_\l)$. Since 
\[\overline{\rho}_{s\otimes \epsilon,2}\iso \overline{\rho}_{s,2}\otimes \epsilon \iso \overline{\rho}_{s,2}\]
is trivial, the section $s\otimes \epsilon$ in fact lifts to $\mathcal{Y}_1(2\l)$ or $\mathcal{Y}(\mu_{2\l})$.
 
Assume that $s$ is rational. Then $s\otimes \epsilon$ is rational by Lemma \ref{lem:beingcuspidaldescends}, hence its lift determines a rational point of $\mathcal{Y}_1(2\l)\iso \mathcal{Y}(\mu_{2\l})$.  By Mazur's Torsion Theorem \cite[Theorem (8)]{mazurModularCurvesEisenstein1977} we have $\mathcal{Y}_1(2\l)(\Q) = \emptyset$, which is a contradiction. Therefore $s$ is cuspidal. 
\ep

\bibliographystyle{amsalpha}

\bibliography{fltforselmersectionsbib}

\end{document}